\documentclass[a4paper,11pt]{amsart}
\usepackage{amsmath,amsthm,amssymb}
\usepackage[mathscr]{eucal}
 \usepackage{cite}
\usepackage{upgreek}
\usepackage[bookmarks,bookmarksnumbered]{hyperref}
\usepackage{enumerate}
\usepackage{amsmath}
\usepackage{mathrsfs}
\usepackage{blindtext}
\usepackage{scrextend}
\usepackage{enumitem}
\addtokomafont{labelinglabel}{\sffamily}
\usepackage{color}
\usepackage[at]{easylist}
\usepackage{bbm}

\numberwithin{equation}{section}

\theoremstyle{plain}
\newtheorem{main}{Theorem}
\newtheorem{mcor}[main]{Corollary}

\newtheorem{theorem}{Theorem}[section]

\newtheorem{lemma}[theorem]{Lemma}
\newtheorem{proposition}[theorem]{Proposition}
\newtheorem{corollary}[theorem]{Corollary}
\theoremstyle{definition}
\newtheorem{definition}[theorem]{Definition}
\newtheorem*{definition*}{Definition}
\newtheorem{example}[theorem]{Example}

\newtheorem{remark}[theorem]{Remark}

\newcommand{\N}{\mathbb{N}}

\newcommand{\C}{\mathbb{C}}

\newcommand{\cA}{\mathcal{A}}

\newcommand{\cR}{\mathcal{R}}
\newcommand{\cS}{\mathcal{S}}
\newcommand{\cT}{\mathcal{T}}
\newcommand{\cU}{\mathcal{U}}

\newcommand{\otb}{\bar{\otimes}}
\newcommand{\dint}{\int^\oplus}
\newcommand{\emb}{\prec}

\newcommand{\Aut}{\operatorname{Aut}}
\newcommand{\Ad}{\operatorname{Ad}}

\newcommand{\eps}{\varepsilon}

\begin{document}

\title[A class of II$_1$ factors with a unique McDuff decomposition]
{A class of II$_1$ factors with a unique McDuff decomposition}

\author[A. Ioana]{Adrian Ioana}
\address{Department of Mathematics, University of California San Diego, 9500 Gilman Drive, La Jolla, CA 92093, USA}
\email{aioana@ucsd.edu}

\author[P. Spaas]{Pieter Spaas}

\address{Department of Mathematics, University of California San Diego, 9500 Gilman Drive, La Jolla, CA 92093, USA}
\email{pspaas@ucsd.edu}

\thanks{The authors were supported in part by NSF Career Grant DMS \#1253402.}
\begin{abstract} 
We provide a fairly large class of II$_1$ factors $N$ such that $M=N\bar{\otimes}R$ has a unique McDuff decomposition, up to isomorphism, where $R$ denotes the hyperfinite II$_1$ factor. This class includes all II$_1$ factors $N=L^{\infty}(X)\rtimes\Gamma$ associated to free ergodic probability measure preserving (p.m.p.) actions $\Gamma\curvearrowright (X,\mu)$ such that either (a) $\Gamma$ is a free group, $\mathbb F_n$, for some $n\geq 2$, or (b) $\Gamma$ is a non-inner amenable group and the orbit equivalence relation of the action $\Gamma\curvearrowright (X,\mu)$ satisfies a property introduced in \cite{JS85}. On the other hand, settling a problem posed by Jones and Schmidt in 1985, we give the first examples of countable ergodic p.m.p. equivalence relations which do not satisfy the property of \cite{JS85}. 
We also prove that if $\mathcal R$ is a countable strongly ergodic p.m.p. equivalence relation and $\mathcal T$ is a hyperfinite ergodic p.m.p. equivalence relation, then $\mathcal R\times\mathcal T$ has a unique stable decomposition, up to isomorphism.
Finally, we provide new characterisations of property Gamma for II$_1$ factors and of strong ergodicity for countable p.m.p. equivalence relations.


\end{abstract}

\maketitle

\section{Introduction and statement of main results}

Let $R$ denote the hyperfinite II$_1$ factor.
The present work is motivated by the following general question: given two II$_1$ factors $N$ and $P$, when are their ``stabilisations" $N\bar{\otimes}R$ and $P\bar{\otimes}R$ isomorphic? 
To put this into context, recall that a II$_1$ factor $M$ is called {\it McDuff} if it absorbs $R$ tensorially, that is, $M\cong M\bar{\otimes}R$. 
This property, introduced by McDuff in \cite{Mc69}, has since played a fundamental role in the theory of von Neumann algebras. In particular, a crucial step in Connes'  celebrated classification of amenable II$_1$ factors \cite{Co75} is showing that any such  factor has McDuff's property.

A II$_1$ factor $M$ is McDuff if and only if it can be decomposed as $M=N\bar{\otimes}R$, for some II$_1$ factor $N$. 
The main result of \cite{Mc69} provides a satisfactory characterisation of when such a decomposition exists: a II$_1$ factor $M$ is McDuff if and only if it admits two non-commuting central sequences.

Our main goal is to investigate the complementary issue of when such a decomposition is unique. 
To make this precise, following \cite{HMV16}, we say 
that a II$_1$ factor $M$ admits a {\it McDuff decomposition} if it can be written as $M=N\bar{\otimes}R$, for some \textit{non-McDuff} II$_1$ factor $N$. 
If two II$_1$ factors $N$ and $P$ are stably isomorphic, that is, if $N\cong P^t$ for some $t>0$, then  $N\bar{\otimes}R$ and $P\bar{\otimes}R$ are isomorphic. 
We say that  a McDuff decomposition $M=N\bar{\otimes}R$ is unique if for any other McDuff decomposition $M=P\bar{\otimes}R$ we necessarily have that $N$ and $P$ are stably isomorphic.
We restrict our attention to II$_1$ factors admitting a McDuff decomposition because if $M$ is a McDuff II$_1$ factor not having a McDuff decomposition (see \cite[Corollary G]{Ma17} for examples), then any II$_1$ factor $N$ satisfying $M\cong N\bar{\otimes}R$ is necessarily isomorphic to $M$ and thus unique, up to isomorphism.

By a striking theorem of Popa \cite[Theorem 5.1]{Po06b}, if  a II$_1$ factor $N$ does not have Murray and von Neumann's {\it property Gamma} \cite{MvN43}, then $M=N\bar{\otimes}R$ admits a unique McDuff decomposition. Moreover, the uniqueness holds up to stable unitary conjugacy: for any other McDuff decomposition $M=P\bar{\otimes}R$, there exists an isomorphism $N\cong P^t$, for some $t>0$, which is implemented by a unitary operator in $M$.  The proof of this result, first presented in \cite{Po04}, relies on Popa's discovery of his influential spectral gap rigidity principle.  In combination with Popa's deformation/rigidity methods this has since led 
to many remarkable applications starting with \cite{Po04,Po06a,Po06b}.

In contrast to \cite{Po06b}, the uniqueness problem for McDuff decompositions as described above where the involved non-McDuff II$_1$ factor $N$ has property Gamma is completely open. 
We make progress on this problem here, by establishing the first unique McDuff decomposition results in the spirit of \cite{Po06b} in the ``property Gamma regime".
Thus, we use methods from Popa's deformation/rigidity theory to provide the first classes of non-McDuff II$_1$ factors $N$ with property Gamma such that  $M=N\bar{\otimes}R$ has a unique McDuff decomposition, see Corollary \ref{B} and Theorem \ref{D}. These classes will be obtained as a consequence of our main technical result:

\begin{main}\label{A}
Let $N$ be a II$_1$ factor which admits a Cartan subalgebra $A$ such that $N'\cap N^{\omega}\subset A^{\omega}$. 
Let $P$ be a non-McDuff II$_1$ factor and $\theta:N\bar{\otimes}R\rightarrow P\bar{\otimes}R$ be an isomorphism.

Then there exist a Cartan subalgebra $B\subset P$, a unitary $u\in P\bar{\otimes}R$, and some $t>0$, such that $P'\cap P^{\omega}\subset B^{\omega}$ and $\theta(A^t)=uBu^*$, where we identify $N\bar{\otimes}R=N^t\bar{\otimes}R^{1/t}$. 

Moreover, $\mathcal R(B\subset P)$ is isomorphic to $\mathcal R(A\subset N)^t$. 
\end{main}

Before giving examples of II$_1$ factors to which Theorem \ref{A} applies, we explain the notions used in its statement.
For the rest of the paper, we fix a free ultrafilter $\omega$  on $\mathbb N$. We denote by  $N^{\omega}$ the {\it ultrapower} of a tracial von Neumann algebra $N$. We say that a von Neumann subalgebra $A$ of $N$ is a {\it Cartan subalgebra} if it is maximal abelian and regular, and denote by $\mathcal R(A\subset N)$ the countable p.m.p. equivalence relation associated to the inclusion $A\subset N$ (see Sections \ref{vN} and \ref{vN2} for details). 

\begin{example} Let $\Gamma$ be a non-inner amenable group and $\Gamma\curvearrowright (X,\mu)$ be a free ergodic p.m.p. action. 
Consider the group measure space von Neumann algebra $N=L^{\infty}(X)\rtimes\Gamma$ associated to the action $\Gamma\curvearrowright (X,\mu)$ \cite{MvN43}. Then $N$ is a II$_1$ factor and $A=L^{\infty}(X)\subset N$ is a Cartan subalgebra.  Moreover, \cite{Ch82} implies that $N'\cap N^{\omega}\subset A^{\omega}$, and thus Theorem \ref{A} applies to $N$. 
If the action $\Gamma\curvearrowright X$ is strongly ergodic, then $N$ does not have property Gamma. In this case, Theorem \ref{A} follows from \cite[Theorem 5.1]{Po06b}, which moreover shows that $\theta(N^t)=uPu^*$, for a unitary $u\in P\bar{\otimes}R$ and some $t>0$. Theorem \ref{A} is new whenever the action $\Gamma\curvearrowright X$ is not strongly ergodic. 

\end{example}


To put Theorem \ref{A} into a better perspective, let $A\subset N$ be an inclusion satisfying its hypothesis, and denote $\mathcal R=\mathcal R(A\subset N)$. A well-known result of Feldman and Moore \cite{FM77} shows that $M$ is isomorphic to the von Neumann algebra $L_{w}(\mathcal R)$ associated to $\mathcal R$ and a $2$-cocycle $w\in$ H$^2(\mathcal R,\mathbb T)$. Theorem \ref{A} thus leads to the following rigidity statement:  any non-McDuff II$_1$ factor $P$ such that $N\bar{\otimes}R\cong P\bar{\otimes}R$ is necessarily isomorphic to $L_v(\mathcal R)^t$, for some $t>0$ and a $2$-cocycle $v\in$ H$^2(\mathcal R,\mathbb T)$. 

While we were unable to determine whether $v$ must be cohomologous to $w$ in general, this is automatically satisfied if the equivalence relation $\mathcal R$ is treeable, leading to the following:

\begin{mcor}\label{B}
Let $n\geq 2$ and $\mathbb F_n\curvearrowright (X,\mu)$ be a free ergodic p.m.p. action. Put $N=L^{\infty}(X)\rtimes\mathbb F_n$. Let $P$ be any II$_1$ factor such that $N\bar{\otimes}R$ and $P\bar{\otimes}R$ are isomorphic.

Then $P$ is either isomorphic to $N^t$, for some $t>0$, or to $N\bar{\otimes}R$.

\end{mcor}

As a particular case of Corollary \ref{B}, we derive a new result for free group measure space factors. Consider a free ergodic p.m.p. action $\mathbb F_p\curvearrowright (Y,\nu)$, for some $p\geq 2$, and $N$ as in the corollary. By a breakthrough theorem of Popa and Vaes \cite{PV11}, $L^{\infty}(X)$ is the unique Cartan subalgebra of $N$, up to unitary conjugacy. Using this result and applying Corollary \ref{B} to $P=L^{\infty}(Y)\rtimes\mathbb F_p$, we deduce that $N\bar{\otimes}R$ and $P\bar{\otimes}R$ are isomorphic if and only if the actions $\mathbb F_n\curvearrowright X$ and $\mathbb F_p\curvearrowright Y$ are stably orbit equivalent. 
A result of Gaboriau \cite{Ga99} further implies that if exactly one of $n$ or $p$ is finite, then $N\bar{\otimes}R\not\cong P\bar{\otimes}R$.

\begin{remark}\label{rk:Ho}
Let $N=L^{\infty}(X)\rtimes\mathbb F_n$, where $\mathbb F_n\curvearrowright (X,\mu)$ is a free ergodic but not strongly ergodic p.m.p. action.  Corollary \ref{B} shows that if we can decompose $N\bar{\otimes}R=P\bar{\otimes}R$, for some non-McDuff II$_1$ factor $P$, then there is an abstract isomorphism between $N^t$ and $P$, for some $t>0$.
This result is optimal in the sense that it cannot be improved to deduce that $N^t$ and $P$ are unitarily conjugate, i.e., that the isomorphism between $N^t$ and $P$ is implemented by a unitary from $N\bar{\otimes}R$. Indeed, since $N$ has property Gamma, \cite[Theorem B]{Ho15} implies that $N\bar{\otimes}R$ admits an (approximately inner) automorphism $\theta$ such that $\theta(N^t)$ is not unitarily conjugate to $N$, inside $N\bar{\otimes}R$, for any $t>0$.
\end{remark}

Before stating the next corollary to Theorem~\ref{A}, we recall that the reason we could not deduce unique McDuff decomposition for every II$_1$ factor to which it applies, is the presence of a 2-cocycle that twists the von Neumann algebra. This difficulty does not appear at the level of the equivalence relations, which allows us to deduce the following:

\begin{mcor}\label{C} 
Let $\mathcal R_1$ and $\mathcal R_2$ be countable ergodic p.m.p. equivalence relations on probability spaces $(X_1,\mu_1)$ and $(X_2,\mu_2)$, respectively. Assume that $L(\mathcal R_1)'\cap L(\mathcal R_1)^{\omega}\subset L^{\infty}(X_1)^{\omega}$ and that $L(\cR_2)$ is not McDuff. Suppose that $\mathcal R_1\times\mathcal T$ is isomorphic to $\mathcal R_2\times\mathcal T$, where $\mathcal T$ is a hyperfinite ergodic p.m.p. equivalence relation on a probability space $(Y,\nu)$.
	
	Then $\mathcal R_2$ is isomorphic to $\mathcal R_1^t$, for some $t>0$.
	
\end{mcor}

Corollary \ref{C} in particular implies that if $\mathcal R_1$ and $\mathcal R_2$ are the orbit equivalence relations of any free ergodic p.m.p. actions of any non-inner amenable groups, then $\mathcal R_1\times\mathcal T$ is isomorphic to $\mathcal R_2\times\mathcal T$ if and only if $\mathcal R_1$ is stably isomorphic to $\mathcal R_2$.

Next, we return to the uniqueness problem for McDuff decompositions of II$_1$ factors. Before providing another class of examples in Theorem~\ref{D} below, we first introduce a property for equivalence relations motivated by a problem posed by Jones and Schmidt. In \cite{JS85} they proved that if $\mathcal S$ is an ergodic but not strongly ergodic (see Section \ref{vN2} for the definition of strong ergodicity) countable p.m.p. 
equivalence relation on a probability space $(X,\mu)$, then $\mathcal S$ admits a hyperfinite quotient. Specifically, there exist a hyperfinite ergodic p.m.p. equivalence relation $\mathcal T$ on a probability space $(Y,\nu)$ and a factor map $\pi:(X,\mu)\rightarrow (Y,\nu)$ such that $\pi([x]_{\mathcal S})=[\pi(x)]_{\mathcal T}$, for almost every $x\in X$  (see \cite[Theorem 2.1]{JS85}).
Jones and Schmidt asked whether one can always find such  $\mathcal T, (Y,\nu), \pi$ with the additional property that $\mathcal S_0=\{(x_1,x_2)\in\mathcal S\mid\pi(x_1)=\pi(x_2)\}$ is strongly ergodic on almost every ergodic component of $\mathcal S_0$ (see \cite[Problem 4.3]{JS85}).
This problem was also mentioned 
 by Schmidt in  \cite[Problem 2.1]{Sc85} as an important question left unanswered by \cite[Theorem 2.1]{JS85}.
Motivated by this problem, we introduce the following:

\begin{definition}\label{JS}
We say that a countable ergodic p.m.p. equivalence relation $\mathcal S$ on a probability space $(X,\mu)$ has the {\it Jones-Schmidt} {\it property} if  there exist a hyperfinite ergodic p.m.p. equivalence relation $\mathcal T$ on a probability space $(Y,\nu)$ and a factor map $\pi:(X,\mu)\rightarrow (Y,\nu)$ such that

\begin{enumerate}
\item $\pi([x]_{\mathcal S})=[\pi(x)]_{\mathcal T}$, for almost every $x\in X$, and
\item the equivalence relation
$\mathcal S_0=\{(x_1,x_2)\in \mathcal S\mid\pi(x_1)=\pi(x_2)\}$ is strongly ergodic on almost every ergodic component of $\mathcal S_0$.
\end{enumerate}

\end{definition}

Using this terminology, \cite[Problem 4.3]{JS85} can be reformulated as follows: does every countable ergodic but not strongly ergodic p.m.p. equivalence relation have the Jones-Schmidt property?

Our next main result shows that $N\otb R$ has a unique McDuff decomposition whenever the above question has a positive answer for the equivalence relation arising from the inclusion of a Cartan subalgebra $A\subset N$ to which Theorem~\ref{A} applies.

\begin{main}\label{D}
Let $N$ be a II$_1$ factor which admits a Cartan subalgebra $A$ such that $N'\cap N^{\omega}\subset A^{\omega}$ and $\mathcal R(A\subset N)$ has the Jones-Schmidt property. 
Let $P$ be any  II$_1$ factor such that $N\bar{\otimes}R$ and $P\bar{\otimes}R$ are isomorphic.

Then $P$ is either isomorphic to $N^t$, for some $t>0$, or to $N\bar{\otimes}R$.

Moreover, assume that $P$ is not McDuff. Then for any isomorphism $\theta:N\bar{\otimes}R\rightarrow P\bar{\otimes}R$, we can find  isomorphisms $\theta_1:N^s\rightarrow P$, $\theta_2:R^{1/s}\rightarrow R$, for some $s>0$, and  an approximately inner automorphism $\Psi: N\bar{\otimes}R\rightarrow N\bar{\otimes}R$ such that $\theta=(\theta_1\otimes\theta_2)\circ\Psi$, where we identify $N\bar{\otimes}R=N^s\bar{\otimes}R^{1/s}$. 
\end{main}

\begin{remark}
The moreover part of Theorem~\ref{D} shows that if $N$ is a II$_1$ factor as in its statement, then $N\otb R$ admits a unique McDuff decomposition, up to ``stable approximately inner conjugacy": given any other McDuff decomposition  $N\otb R=P\otb R$, there exists an isomorphism $N^s\cong P$, for some $s>0$, which  is implemented by an approximately inner automorphism of $N\bar{\otimes}R$. As in Remark~\ref{rk:Ho}, by \cite[Theorem~B]{Ho15} this is optimal and cannot be improved to an implementation by an inner automorphism.
\end{remark}

\begin{example} Theorem \ref{D} applies to a large class of group measure space II$_1$ factors. To see this,
let $\Gamma$ be an infinite group,  $\Sigma$ an infinite amenable group, $\delta:\Gamma\rightarrow\Sigma$ an onto group homomorphism, and $\Gamma\curvearrowright (Z,\eta)$, $\Sigma\curvearrowright (Y,\nu)$ be free ergodic p.m.p. actions such that the action $\ker{\delta}\curvearrowright (Z,\eta)$ is strongly ergodic. Let $(X,\mu)=(Y,\nu)\times (Z,\eta)$ and consider the action $\Gamma\curvearrowright (X,\mu)$ given by $g\cdot (y,z)=(\delta(g)y,gz)$, for all $g\in\Gamma, y\in Y, z\in Z$.

 Then  $\mathcal S=\mathcal R(\Gamma\curvearrowright X):=\{(x_1,x_2)\in X\times X\mid\Gamma\cdot x_1=\Gamma\cdot x_2\}$ has the Jones-Schmidt property. To see this, let $\mathcal T=\mathcal R(\Sigma\curvearrowright Y)$, 
and define the factor map $\pi:X\rightarrow Y$ by $\pi(y,z)=y$. 
Since $\Sigma$ is infinite amenable and the action $\Sigma\curvearrowright Y$ is ergodic, $\mathcal T$ is a hyperfinite ergodic p.m.p. equivalence relation by \cite{OW80}. 
Also, $\pi([x]_{\mathcal S})=[\pi(x)]_{\mathcal T}$ for all $x\in X$.
Moreover, if $\mathcal S_0=\{(x_1,x_2)\in\mathcal S\mid\pi(x_1)=\pi(x_2)\}$, then $\mathcal S_0$ is equal to $\{(y,y)|y\in Y\}\times \mathcal R(\ker\delta\curvearrowright Z)$.
Thus, every ergodic component of $\mathcal S_0$ is isomorphic to $\mathcal R(\ker\delta\curvearrowright Z)$, and hence is strongly ergodic.

Moreover, assume that $\Gamma$ is not inner amenable, denote $N=L^{\infty}(X)\rtimes\Gamma$ and $A=L^{\infty}(X)$. Then $N'\cap N^{\omega}\subset A^{\omega}$ by \cite{Ch82}, and therefore Theorem \ref{D} applies to $N$.
\end{example}

The above example gives a wide class of examples of countable ergodic p.m.p. equivalence relations with the Jones-Schmidt property. We next show that there exist countable ergodic but not strongly ergodic p.m.p. equivalence relations that fail to have the Jones-Schmidt property. 
 This settles the question asked by Jones and Schmidt in \cite[Problem~4.3]{JS85}. 
More specifically, in Theorems~\ref{E} and \ref{E'} below, we give two large classes of free ergodic p.m.p. actions whose orbit equivalence relations do not have the Jones-Schmidt property. First, we prove that this holds for any infinite product of strongly ergodic actions:
\begin{main}\label{E}
For every $n\in\mathbb N$, let $\Gamma_n\curvearrowright (X_n,\mu_n)$ be a strongly ergodic free p.m.p. action of an infinite countable group $\Gamma_n$. Define $\Gamma=\bigoplus_{n\in\mathbb N}\Gamma_n$, $(X,\mu)=\prod_{n\in\mathbb N}(X_n,\mu_n)$, and consider the product action $\Gamma\curvearrowright (X,\mu)$ given by $g\cdot x=(g_n\cdot x_n)_n$, for all $g=(g_n)_n\in\Gamma$ and $x=(x_n)_n\in X$.

Then the orbit equivalence relation $\mathcal R(\Gamma\curvearrowright X)$ does not have the Jones-Schmidt property.
\end{main}
 
 A countable group $\Lambda$ admits strongly ergodic free p.m.p. actions if and only if it is non-amenable. Moreover, if $\Lambda$ is non-amenable, then the Bernoulli action $\Lambda\curvearrowright (Y_0,\nu_0)^{\Lambda}$ is strongly ergodic, for any non-trivial probability space $(Y_0,\nu_0)$.
Theorem \ref{E} therefore applies to any infinite direct sum $\Gamma=\bigoplus_{n\in\mathbb N}\Gamma_n$ of non-amenable groups and shows that any such group admits free ergodic p.m.p. actions whose orbit equivalence relations fail the Jones-Schmidt property.
Next, we provide a second such result which, unlike Theorem \ref{E}, covers actions of non-abelian free groups (see Example \ref{free}).

\begin{main}\label{E'} Let $\Gamma$ be a countable group. For every $n\in\mathbb N$, let $\Gamma\curvearrowright (X_n,\mu_n)$ be a free ergodic p.m.p. action such  that  the diagonal action $\Gamma\curvearrowright (X_n\times X_n,\mu_n\times\mu_n)$ has spectral gap. Assume that we can find $F_{n,k}\in L^{\infty}(X_n)$, for all $n,k\in\mathbb N$, such that  $\sup_{n,k}\|F_{n,k}\|_{\infty}\leq 1$, $\inf_{n,k}\|F_{n,k}\|_2>0$,

\begin{itemize}\item $F_{n,k}\rightarrow 0$ weakly in $L^2(X_n)$, as $k\rightarrow\infty$, for every $n\in\mathbb N$, and \item $\lim\limits_{n\rightarrow\infty}\big(\sup_{k\in\mathbb N}\|F_{n,k}\circ g-F_{n,k}\|_2\big)=0$, for every $g\in\Gamma$.\end{itemize}
Consider the diagonal action $\Gamma\curvearrowright (X,\mu)=\prod_{n\in\mathbb N}(X_n,\mu_n)$ given by $g\cdot x=(g\cdot x_n)_n$, for all $g\in\Gamma$ and $x=(x_n)_n\in X$.

Then the orbit equivalence relation $\mathcal R(\Gamma\curvearrowright X)$ does not have the Jones-Schmidt property.
\end{main}

\begin{example}\label{free} Let $\Gamma=\mathbb F_m$ be the free group on $2\leq m\leq\infty$ generators.
Denote by $|g|$ the word length of an element $g\in\Gamma$ with respect to a free set of generators. Let $t>0$. By a well-known result of Haagerup \cite{Ha79},  the function $\varphi_t:\Gamma\rightarrow\mathbb R$ given by $\varphi_t(g)=e^{-t|g|}$ is positive definite. Denote by $\pi_t:\Gamma\rightarrow\mathcal O(\mathcal H_t)$ the GNS orthogonal representation associated to $\varphi_t$. Let $\tilde{\mathcal H}_t=\oplus_{i\in\mathbb N}\mathcal H_t$ and  $\tilde{\pi}_t=\oplus_{i\in\mathbb N}\pi_t:\Gamma\rightarrow\mathcal O(\tilde{\mathcal H}_t)$ be the direct sum of infinitely many copies of $\pi_t$. 
Let $\Gamma\curvearrowright (X_t,\mu_t)$ be the Gaussian action associated to $\widetilde{\pi}_t$  (see, e.g., \cite[Chapter 3]{Gl03})
Let $(t_n)$ be any sequence of positive numbers converging to $0$. As we will prove in Remark \ref{freeproof}, the diagonal action of $\Gamma$ on $(X,\mu)=\prod_{n\in\mathbb N}(X_{t_n},\mu_{t_n})$ satisfies 
the hypothesis of Theorem \ref{E'}.
\end{example}

Let $\cT$ denote a hyperfinite ergodic p.m.p. equivalence relation. Similar to the von Neumann algebra case considered above, one can ask the following question: given two countable ergodic p.m.p. equivalence relations $\cR_1$ and $\cR_2$, when are their ``stabilisations'' $\cR_1\times \cT$ and $\cR_2\times \cT$ isomorphic? We recall from \cite{JS85} that a countable ergodic p.m.p. equivalence relation $\cS$ is called \textit{stable} if it can be decomposed as $\cS=\cR\times \cT$, for some countable ergodic equivalence relation $\cR$. In \cite[Theorem 3.4]{JS85} stability is shown to be equivalent to an asymptotic property: 
$\mathcal S$ is stable if and only if its full group contains a nontrivial asymptotically central sequence. In recent years, there has been a surge of interest in the study of stability for equivalence relations, see, e.g., the articles \cite{Ki12,TD14,Ki16,DV16,Ma17}. 

We contribute to this study here by investigating the problem of when a decomposition as above  is unique. By analogy with the case of McDuff II$_1$ factors, we say that a countable ergodic p.m.p. equivalence relation $\cS$ admits a \textit{stable decomposition} if it can be written as $\cS=\cR\times\cT$ for some non-stable countable ergodic p.m.p. equivalence relation $\cR$. If two countable ergodic p.m.p. equivalence relations $\cR_1$ and $\cR_2$ are stably isomorphic, that is, if $\cR_1\cong \cR_2^t$ for some $t>0$, 
then $\cR_1\times\cT$ and $\cR_2\times\cT$ are isomorphic. We say that a stable decomposition $\cS=\cR_1\times \cT$ is unique if for any other stable decomposition $\cS=\cR_2\times \cT$ we necessarily have that $\cR_1$ and $\cR_2$ are stably isomorphic.

Note that Corollary~\ref{C} already gave the first result in this direction: let $\cR_1$ be a countable ergodic p.m.p. equivalence relation on a probability space $(X_1,\mu_1)$ such that $L(\mathcal R_1)'\cap L(\mathcal R_1)^{\omega}\subset L^{\infty}(X_1)^{\omega}$. Then any countable ergodic p.m.p. equivalence relation $\mathcal R_2$ whose von Neumann algebra $L(\mathcal R_2)$ is not McDuff and satisfies $\cR_1\times \cT\cong \cR_2\times \cT$ must be stably isomorphic to $\mathcal R_1$. 
However, by \cite{CJ82}, there exist equivalence relations $\mathcal R_2$ which are strongly ergodic, and thus not stable, such that $L(\mathcal R_2)$ is McDuff. As such, Corollary \ref{C} does not imply that $\mathcal R_1\times\mathcal T$ has a unique stable decomposition.

Our next main result completely settles the uniqueness problem for stable decompositions of $\mathcal R_1\times\mathcal T$ under the assumption that $\cR_1$ is strongly ergodic. More precisely, we have the following:

\begin{main}\label{F}
	Let $\cR_1$ and  $\cR_2$ be countable ergodic p.m.p. equivalence relations on probability spaces $(X_1,\mu_1)$ and $(X_2,\mu_2)$, respectively. Assume that $\cR_1$ is strongly ergodic.
	Suppose that $\mathcal R_1\times\mathcal T$ is isomorphic to $\mathcal R_2\times\mathcal T$, where $\mathcal T$ is a hyperfinite ergodic p.m.p. equivalence relation on a probability space $(Y,\nu)$.

	Then either
	\begin{enumerate}
		\item $\cR_2$ is also strongly ergodic and $\cR_2\cong \cR_1^t$, for some $t>0$, or
		\item$\cR_2$ is stable and $\cR_2\cong \cR_1\times \cT$. 
	\end{enumerate}
\end{main}

Finally, we provide new characterisations of property Gamma for II$_1$ factors and of strong ergodicity for countable ergodic p.m.p. equivalence relations. We include these results here,  although they are not related to the results stated above, as they appear to be of independent interest.

On \cite[page 92]{JS85}, Jones and Schmidt pointed out that their characterisation of strong ergodicity for countable equivalence relations \cite[Theorem 2.1]{JS85} has no obvious analogue in the setting of von Neumann algebras. We  address this problem here by giving 
such an analogue. In particular, we show that  any II$_1$ factor $M$ with property Gamma admits a regular diffuse abelian von Neumann subalgebra.

\begin{main}\label{G}
Let $M$ be a separable II$_1$ factor. Then the following are equivalent:
\begin{enumerate}
\item $M$ has property Gamma.
\item There exist a hyperfinite subfactor $R\subset M$ and a Cartan subalgebra $A\subset R$ such that $M$ is generated by $R$ and $A'\cap M$. 
\item There exists a regular abelian von Neumann subalgebra $A\subset M$ such that $\mathcal R(A\subset M)$ is hyperfinite and ergodic. 
\end{enumerate}
\end{main}

In order to state our last result, we need some additional terminology.
Let $\mathcal S$ be a countable p.m.p. equivalence relation on a probability space $(X,\mu)$ and $E$ be a Borel equivalence relation on a standard Borel space $Y$. A map $\varphi:X\rightarrow Y$ is a {\it homomorphism} from $\mathcal S$ to $E$ if $(\varphi\times\varphi)(\mathcal S)\subset E$.
We say that $\mathcal S$ is {\it $E$-ergodic} if  for any Borel homomorphism $\varphi:X\rightarrow Y$, there exists $y\in Y$ such that $(\varphi(x),y)\in E$, for almost every $x\in X$. 
Let $E_0$ be the hyperfinite Borel equivalence relation on $\{0,1\}^{\mathbb N}$ given by $(y_n)E_0(z_n)$ if and only if $y_n=z_n$, for all but finitely many $n$. As shown in \cite[Theorem A.2.2]{HK05}, a countable ergodic p.m.p. equivalence relation $\mathcal S$ is strongly ergodic if and only if it is $E_0$-ergodic. We prove that strong ergodicity can also be characterized as $E$-ergodicity, where $E$ is the orbit equivalence relation of the left translation action $\text{Inn}(R)\curvearrowright \text{Aut}(R)$. 


\begin{main}\label{H} Let $\mathcal S$ be a countable ergodic p.m.p. equivalence relation on a probability space $(X,\mu)$. 
Let $\mathcal T$ be a countable hyperfinite ergodic p.m.p. equivalence relation on a probability space $(Y,\nu)$. Let $R$ denote the hyperfinite II$_1$ factor.
Then the following are equivalent:
\begin{enumerate}
\item $\mathcal S$ is strongly ergodic.
\item For any measurable map $\varphi:X\rightarrow\text{Aut}(\mathcal T)$ satisfying $\varphi(x_1)^{-1}\varphi(x_2)\in [\mathcal T]$ for almost every $(x_1,x_2)\in\mathcal S$, there exists $\alpha\in\text{Aut}(\mathcal T)$ such that $\varphi(x)\alpha\in [\mathcal T]$, for almost every $x\in X$.
\item For any measurable map $\varphi:X\rightarrow\text{Aut}(R)$ satisfying $\varphi(x_1)^{-1}\varphi(x_2)\in\text{Inn}(R)$ for almost every $(x_1,x_2)\in\mathcal S$, there exists $\alpha\in\text{Aut}(R)$ such that $\varphi(x)\alpha\in\text{Inn}(R)$, for almost every $x\in X$.
\end{enumerate}

\end{main}

\subsection*{Organization of the paper} Besides the introduction this paper has seven other sections. First of all, in the Preliminaries, we recall some tools needed in the remainder of the paper and prove some useful lemmas. In Section~3 we prove the main technical result, Theorem~\ref{A}, and derive Corollaries~\ref{B} and \ref{C} from it. 
We then continue in Section~4 with proving a cocycle rigidity result which, besides being of independent interest, will be used in the proofs of Theorems~\ref{D} and \ref{H}. In Section~5 we investigate the Jones-Schmidt property. We provide an (a priori stronger) equivalent characterisation and use this to prove Theorem~\ref{D}. Sections~6 and ~7 are devoted to the proofs of Theorems~\ref{E}, \ref{E'} and \ref{F}. 
Finally, in Section~8 we prove 
Theorems~\ref{G} and \ref{H}.

\subsection*{Acknowledgements} The first author would like to thank Lewis Bowen and Vaughan Jones for stimulating discussion regarding \cite[Problem 4.3]{JS85}. In particular, he is grateful to Vaughan for pointing out the analogy between this problem and a problem of Connes solved by Popa in \cite{Po09}, and to Lewis for raising a question which motivated Theorem \ref{E'}.  We would also like to thank Sorin Popa for several comments on the first version of this paper.
Part of this work was done during the program ``Quantitative Linear Algebra" at the Institute for Pure and Applied Mathematics. We would like to thank IPAM for its hospitality.

\section{Preliminaries} We begin by fixing notation regarding von Neumann algebras and equivalence relations.

\subsection{von Neumann algebras}\label{vN}  In this paper, we consider {\it tracial von Neumann algebras} $(M,\uptau)$, i.e., von Neumann algebras $M$ endowed with a faithful normal tracial state $\uptau:M\rightarrow\mathbb C$. For $x\in M$, we denote by $\|x\|$ its operator norm and by $\|x\|_2:=\uptau(x^*x)^{1/2}$ its (so-called) $2$-norm. We denote  by $(M)_1$ the set of $x\in M$ with $\|x\|\leq 1$, and by $L^2(M)$ the Hilbert space obtained by taking the closure of $M$ with respect to the $2$-norm. 
Unless stated otherwise, we will always assume that $M$ is a {\it separable} von Neumann algebra, i.e.,
 $L^2(M)$ is a separable Hilbert space.
We denote by $\mathcal U(M)$ the group of unitaries of $M$, by Aut$(M)$ the group of $\uptau$-preserving automorphisms of $M$, and by Inn$(M)$ the group of {\it inner} automorphisms of $M$, i.e. those of the form $\text{Ad}(u)(x)=uxu^*$, where $u\in\mathcal U(M)$. 
We endow Aut$(M)$ with the Polish topology given by  pointwise $\|.\|_2$-convergence. If $(X,\mu)$ is a standard measure space and $G$ is a Polish group (e.g., $\mathcal U(M)$ or Aut$(M)$), we say that a map $f:X\rightarrow G$ is measurable if $f^{-1}(B)$ is a measurable subset of $X$, for every Borel subset $B\subset G$.
For a subgroup $H<\text{Aut}(M)$, we denote $M^H=\{x\in M\mid\text{$\theta(x)=x$, for all $\theta\in H$}\}$.

Let $P,Q\subset M$ be von Neumann subalgebras, which we will always assume to be unital. We denote by  $E_P:M\rightarrow P$ the unique $\uptau$-preserving {\it conditional expectation} from $M$ onto $P$. We write that $Q\subset_{\varepsilon}P$, for some $\varepsilon>0$, if we have  $\|x-E_P(x)\|_2\leq\varepsilon$, for every $x\in (Q)_1$. 
 We denote by $P'\cap M=\{x\in M\mid\text{$xy=yx$, for all $y\in P$}\}$ the {\it relative commutant} of $P$ in $M$, and by $\mathcal N_M(P)=\{u\in\mathcal U(M)\mid uPu^*=P\}$ the {\it normalizer} of $P$ in $M$. We say that $P$ is {\it regular} in $M$ if the von Neumann algebra generated by $\mathcal N_M(P)$ is equal to $M$.  We denote by $\mathcal Z(M)=M'\cap M$ the {\it center} of $M$, and by $P\vee Q$ the  von Neumann subalgebra of $M$ generated by $P$ and $Q$.

Let $\omega$ be a free ultrafilter on $\mathbb N$. Consider the C$^*$-algebra $\ell^{\infty}(\mathbb N,M)=\{(x_n)\in M^{\mathbb N}\mid\sup\|x_n\|<\infty\}$ together with its closed ideal $$\mathcal I_{\omega}=\{(x_n)\in\ell^{\infty}(\mathbb N,M)\mid\lim\limits_{n\rightarrow\omega}\|x_n\|_2=0\}.$$  Then the C$^*$-algebra $M^{\omega}:=\ell^{\infty}(\mathbb N,M)/\mathcal I_{\omega}$ is in fact a tracial von Neumann algebra, called the {\it ultrapower } of $M$, whose trace is given by $\uptau_{\omega}(x)=\lim\limits_{n\rightarrow\omega}\uptau(x_n)$, for all $x=(x_n)\in M^{\omega}$. Note that if $M$ is diffuse (or has a diffuse direct summand), then $M^{\omega}$ is non-separable. For an automorphism $\theta$ of $M$ we denote still by $\theta$ the automorphism of $M^{\omega}$ given by $\theta((x_n))=(\theta(x_n))$. In this manner, we view any subgroup $G<\text{Aut}(M)$ as a subgroup $G<\text{Aut}(M^{\omega})$.
Finally, if $M_n, n\in\mathbb N$, is a sequence of von Neumann subalgebras of $M$, then their {\it ultraproduct}, denoted by $\prod_{\omega}M_n$, can be realized  as the von Neumann subalgebra of $M^{\omega}$ consisting of $x=(x_n)$ such that $\lim\limits_{n\rightarrow\omega}\|x_n-E_{M_n}(x_n)\|_2=0$.

\subsection{Equivalence relations and Cartan subalgebras}\label{vN2}
Let $(X,\mu)$ be a probability space, which we will always assume to be {\it standard}. 
If $\Gamma\curvearrowright (X,\mu)$ is a p.m.p. action of a countable group $\Gamma$, then its {\it orbit equivalence relation} $\mathcal R(\Gamma\curvearrowright X)=\{(x,y)\in X^2\mid\Gamma\cdot x=\Gamma\cdot y\}$ is countable p.m.p. Conversely, every countable p.m.p. equivalence relation arises in this way by \cite{FM77}.

Let $\mathcal R$ be a countable p.m.p. equivalence relation on $(X,\mu)$. For every $x\in X$, we denote by $[x]_{\mathcal R}$ its equivalence class. We endow $\mathcal R$ with an infinite measure $\bar{\mu}$ given by $$\text{$\bar{\mu}(A)=\int_{X}\#\{y\in X\mid (x,y)\in A\}\;\text{d}\mu(x)$, for every Borel set $A\subset\mathcal R$.}$$
 The
{\it automorphism group} of $\mathcal R$, denoted by Aut$(\mathcal R)$, consists of all measure space automorphisms $\alpha$ of $(X,\mu)$ such that $(\alpha(x),\alpha(y))\in\mathcal R$, for $\bar{\mu}$-almost every $(x,y)\in\mathcal R$.
The {\it full group} of $\mathcal R$, denoted by $[\mathcal R]$, is the subgroup of  all $\alpha\in$ Aut$(\mathcal R)$  such that $(\alpha(x),x)\in\mathcal R$, for $\mu$-almost every $x\in X$. The {\it full pseudogroup} of $\mathcal R$, denoted by $[[\mathcal R]]$, consists of all isomorphisms $\alpha:(A,\mu_{|A})\rightarrow (B,\mu_{|B})$, where $\mu_A,\mu_B$ denote the restrictions of $\mu$ to  measurable sets $A, B\subset X$, such that $(\alpha(x),x)\in \mathcal R$, for almost every $x\in A$. Given a Polish group $G$, a measurable map $c:\mathcal R\rightarrow G$ is called  a $1$-{\it cocycle} if it satisfies the identity $c(x,y)=c(x,z)c(z,y)$, for all $y,z\in [x]_{\mathcal R}$, for $\mu$-almost every $x\in X$. 

We say that $\mathcal R$ is {\it ergodic} if for every measurable set $A\subset X$ satisfying  $\mu(\alpha(A)\triangle A)=0$, for all $\alpha\in [\mathcal R]$, we have $\mu(A)\in\{0,1\}$. We say that $\mathcal R$ is {\it strongly ergodic} if for every sequence of measurable sets $A_n\subset X$ satisfying $\lim\limits_{n\rightarrow\infty}\mu(\alpha(A_n)\triangle A_n)=0$, for all $\alpha\in [\mathcal R]$, we have $\lim\limits_{n\rightarrow\infty}\mu(A_n)(1-\mu(A_n))=0$. Note that a p.m.p. action $\Gamma\curvearrowright (X,\mu)$ of a countable group $\Gamma$ is ergodic (or, strongly ergodic) if and only if its orbit equivalence relation $\mathcal R(\Gamma\curvearrowright X)$ is.

Assume that $\mathcal R$ is ergodic, and let $t>0$. Denote by $\#$ the counting measure of $\mathbb N$, and consider a measurable set $X^t\subset X\times\mathbb N$ with $(\mu\times\#)(X^t)=t$. The $t$-{\it amplification} of $\mathcal R$, denoted by $\mathcal R^t$, is defined as  the equivalence relation on $X^t$ given by $((x,m),(y,n))\in\mathcal R^t\Longleftrightarrow (x,y)\in\mathcal R$. Since $\mathcal R$ is ergodic, the isomorphism class of $\mathcal R^t$ only depends on $\mathcal R$ and $t$, but not on the choice of $X^t$.

Let $(M,\uptau)$ be a separable tracial von Neumann algebra and $A\subset M$ be an abelian von Neumann subalgebra. Identify $A=L^{\infty}(X)$, for some standard probability space $(X,\mu)$. Then for every $u\in\mathcal N_{M}(A)$, we can find an automorphism $\alpha_u$ of $(X,\mu)$ such that $a\circ\alpha_u=uau^*$, for every $a\in A$. The equivalence relation $\mathcal R(A\subset M)$ of the inclusion $A\subset M$ is the smallest countable p.m.p. equivalence relation on $(X,\mu)$ whose full group contains $\alpha_u$, for every $u\in\mathcal N_M(A)$.

Now, assume that $M$ is a II$_1$ factor and $A\subset M$ is a {\it Cartan subalgebra}, i.e. a maximal abelian regular von Neumann subalgebra. Then $\mathcal R:=\mathcal R(A\subset M)$ is ergodic and there exists a 2-cocycle $w\in$ H$^2(\mathcal R,\mathbb T)$ such that the inclusions $(A\subset M)$ and $(L^{\infty}(X)\subset L_{w}(\mathcal R))$ are isomorphic \cite{FM77}. For $t>0$, 
the $t$-{\it amplification} of $M$, denoted by $M^t$, is defined as the isomorphism class of $p(\mathbb B(\ell^2)\bar{\otimes}M)p$, where $p\in\mathbb B(\ell^2)\bar{\otimes}M$ is a projection with $(\text{Tr}\otimes\uptau)(p)=t$, $\text{Tr}$ denoting the usual trace on $\mathbb B(\ell^2)$.
Similarly, the inclusion $(A^t\subset M^t)$ is defined as the isomorphism class of the inclusion $(p(\ell^{\infty}\otimes A)p\subset p(\mathbb B(\ell^2)\bar{\otimes}M)p)$, where  $p\in\mathbb B(\ell^2)\bar{\otimes}A$ is a projection with $(\text{Tr}\otimes\uptau)(p)=t$, and $\ell^{\infty}\subset\mathbb B(\ell^2)$ is the subalgebra of diagonal operators.
With this terminology, we have that $A^t\subset M^t$ is a Cartan subalgebra, and $\mathcal R(A^t\subset M^t)\cong\mathcal R(A\subset M)^t$.

\subsection{Popa's intertwining-by-bimodules} 
In \cite{Po03}, Popa introduced a powerful theory for deducing unitary conjugacy of subalgebras of tracial von Neumann algebras. This theory, which we recall next, will be a key tool in the present paper. 

Let $P,Q$ be  von Neumann subalgebras of a tracial von Neumann algebra $(M,\uptau)$.  We say that {\it $P$ embeds into $Q$ inside $M$} and write $P\prec_{M}Q$ if we can find non-zero projections $p\in P,q\in Q$, a $*$-homomorphism $\theta:pPp\rightarrow qQq$, and a non-zero partial isometry $v\in qMp$ satisfying $\theta(x)v=vx$, for all $x\in pPp$. 
Moreover, if $Pp'\prec_{M}Q$, for any non-zero projection $p'\in P'\cap M$, we write $P\prec_{M}^{s}Q$.

\begin{theorem} [\!\!\cite{Po03}]\label{corner}
In the above setting, let $\mathcal U\subset\mathcal U(P)$ be a subgroup which generates $P$. Then the following conditions are equivalent:
\begin{enumerate}
\item $P\prec_MQ$.
\item There is no sequence $u_n\in \mathcal U$ satisfying $\|E_Q(au_nb)\|_2\rightarrow 0$, for all $a,b\in M$.

\end{enumerate} 
\end{theorem}

If $P\subset_{\varepsilon}Q$, for some $\varepsilon<1$, then Theorem \ref{corner} implies that  $P\prec_{M}Q$. Moreover, this can be made quantitative, as follows:

\begin{lemma}\label{epsilon}
In the above setting, assume that $P\subset_{\varepsilon}Q$, for some $\varepsilon>0$.
Then there exists a projection $p'\in\mathcal Z(P'\cap M)$ such that $Pp'\prec_{M}^{s}Q$ and $\uptau(p')\geq 1-\varepsilon.$ 
\end{lemma}

{\it Proof.} Let $p'\in\mathcal Z(P'\cap M)$ be the maximal projection such that $Pp'\prec_{M}^s Q$. If we let $p''=1-p'$, then by \cite[Lemma 2.3(3)]{DHI16} we get that $Pp''\nprec_{M}Q$. On the other hand, for every $u\in\mathcal U(P)$ we have that $\|u-E_Q(u)\|_2\leq\varepsilon$, hence $\Re\uptau(up''E_Q(u)^*)=\uptau(p'')+\Re\uptau(up''(E_Q(u)-u)^*)\geq\uptau(p'')-\varepsilon$, and therefore $\|E_Q(up'')\|_2\geq\uptau(p'')-\varepsilon$, for all $u\in\mathcal U(P)$.
If $\uptau(p'')-\varepsilon>0$, then by Theorem \ref{corner} we would get that $Pp''\prec_{M}Q$, which is a contradiction. Hence, $\uptau(p'')\leq\varepsilon$ and thus $\uptau(p')\geq 1-\varepsilon$. \hfill$\blacksquare$

\subsection{Relative commutants in ultrapower II$_1$ factors} In this subsection, we establish several results concerning relative commutants in ultrapower II$_1$ factors. First, we record the following  easy and well-known fact that will be used repeatedly.

\begin{lemma}\label{easy}
Let $(M,\uptau)$ be a tracial von Neumann algebra, and  $P\subset M$, $B_n\subset M$, for $n\in\mathbb N$, be von Neumann subalgebras.
Assume that $\prod_{\omega}B_n\subset P^{\omega}$. Then there exists a sequence of positive real numbers $(\varepsilon_n)$ such that $B_n\subset_{\varepsilon_n}P$, for every $n\in\mathbb N$, and  $\lim\limits_{n\rightarrow\omega}\varepsilon_n=0$. 
\end{lemma}

{\it Proof.} For every $n\in\mathbb N$, let $\varepsilon_n=\sup\{\|x-E_P(x)\|_2|\;x\in (B_n)_1\}$, and $x_n\in (B_n)_1$ such that $\|x_n-E_{P}(x_n)\|\geq\varepsilon_n-2^{-n}$. Let $x=(x_n)\in\prod_{\omega}B_n$. Then $\lim\limits_{n\rightarrow\omega}\|x_n-E_P(x_n)\|_2=\|x-E_{P^{\omega}}(x)\|_2=0$, which implies that $\lim\limits_{n\rightarrow\omega}\varepsilon_n=0$.
\hfill$\blacksquare$

The next result is a particular case of \cite[Lemma 2.2]{CD18}.
 For completeness, we include a short argument based on Lemma \ref{easy}.

\begin{lemma}[\!\!\cite{CD18}]\label{intert}
Let $(M,\uptau)$ be a tracial von Neumann algebra, $R\subset M$ be the hyperfinite II$_1$ factor, and $P\subset M$ a von Neumann subalgebra.
Assume that $R'\cap R^{\omega}\subset P^{\omega}$.
Then $R\prec_{M}^{s}P$.
\end{lemma}

{\it Proof.} Let $p\in R'\cap M$ be a non-zero projection.  Write $R=\bar{\otimes}_{k\in\mathbb N}\mathbb M_2(\mathbb C)$, and put $R_n=\bar{\otimes}_{k\geq n}\mathbb M_2(\mathbb C)$, for $n\in\mathbb N$. Since $\prod_{\omega}R_n\subset R'\cap R^{\omega}$,  Lemma \ref{easy} provides $n_0\geq 1$ such that $\|u-E_{P}(u)\|_2\leq \uptau(p)/2$, for all $u\in\mathcal U(R_{n_0})$. This implies that $\|E_P(up)\|_2\geq\uptau(upE_P(u)^*)\geq\uptau(p)/2$, for all $u\in\mathcal U(R_{n_0})$. 
Thus, $R_{n_0}p\prec_{M}P$, and since $R_{n_0}p\subset Rp$ is a finite index subfactor, we get that $Rp\prec_{M}P$. 
\hfill$\blacksquare$

The following result is a particular case of Ocneanu's central freedom lemma, see \cite[Lemma 15.25]{EK97}.
Nevertheless, we include a proof for the reader's convenience. 

\begin{theorem}[\!\!\cite{EK97}]\label{comm}
Let $(M,\uptau)$ be a tracial von Neumann algebra and $R\subset M$ be the hyperfinite II$_1$ factor.
Then $(R'\cap R^{\omega})'\cap M^{\omega}=(R'\cap M)^{\omega}\vee R$.
\end{theorem}

{\it Proof.} Since clearly $(R'\cap M)^{\omega}\vee R\subset (R'\cap R^{\omega})'\cap M^{\omega}$, we only have to show the reverse inclusion. To this end, 
we adapt the proof of \cite[Theorem 2.1.2]{Po13} which deals with the case $M=R$. 
  Write $R=\bar{\otimes}_{k\in\mathbb N}\mathbb M_2(\mathbb C)$. Define $R_m=\bar{\otimes}_{k\geq m}\mathbb M_2(\mathbb C)$ and $R^m=\bar{\otimes}_{k<m}\mathbb M_2(\mathbb C)$, for every $m\in\mathbb N$.

{\bf Claim}.
For all $x\in (R'\cap R^{\omega})'\cap M^{\omega}$ and $\varepsilon>0$, there exists $k\in\mathbb N$ such that $\|x-E_{(R_k'\cap M)^{\omega}}(x)\|_2\leq\varepsilon$. 

{\it Proof of the claim.} We may assume that $\lim\limits_{n\rightarrow\omega}\|x_n-E_{R_1'\cap M}(x_n)\|_2=\|x-E_{(R_1'\cap M)^{\omega}}(x)\|_2>\varepsilon$, otherwise the claim holds for $k=1$. If we let $V=\{n\in\mathbb N|\;\|x_n-E_{R_1'\cap M}(x_n)\|_2>\varepsilon\}$, then $V\in\omega$.
Fix $n\in V$, and denote $F_n=\{k\in\mathbb N\mid \|x_n-E_{R_k'\cap M}(x_n)\|_2>\varepsilon\}$. Then $1\in F_n$, and we define $$k(n)=\begin{cases}\max F_n\;\;\;\text{if $F_n$ is finite}\\ n\;\;\;\text{if $F_n$ is infinite.}
\end{cases}$$

We claim that there is $u_n\in\mathcal U(R_{k(n)})$ such that $\|x_n-u_nx_nu_n^*\|_2\geq\varepsilon/2$. Suppose  by contradiction that $\|x_n-ux_nu^*\|_2<\varepsilon/2$, for all $u\in\mathcal U(R_{k(n)})$. Let $C_n\subset L^2(M)$ be the $\|.\|_2$-closure of the convex hull of $\{ux_nu^*| u\in\mathcal U(R_{k(n)})\}$. Then $\|x_n-\xi\|_2\leq\varepsilon/2$, $\|\xi\|\leq \|x_n\|$, and $E_{R_{k(n)}'\cap M}(\xi)=E_{R_{k(n)}'\cap M}(x_n)$, for all $\xi\in C_n$.
Let $\eta$ be the unique element of minimal $\|.\|_2$ in $C_n$. Since $u\xi u^*\in C_n$ and $\|u\xi u^*\|_2=\|\xi\|_2$, for all $\xi\in C_n$ and $u\in\mathcal U(R_{k(n)})$, the uniqueness of $\eta$ implies that $\eta=u\eta u^*$, for all $u\in\mathcal U(R_{k(n)})$. Thus, $\eta\in R_{k(n)}'\cap M$ and hence $\eta=E_{R_{k(n)}'\cap M}(\eta)=E_{R_{k(n)}'\cap M}(x_n)$. This however contradicts the fact that $\|x_n-E_{R_{k(n)}'\cap M}(x_n)\|_2\geq\varepsilon$, by the construction of $k(n)$.

For  every $n\in\mathbb N\setminus V$, we put $k(n)=1$ and $u_n=1$. We claim that $\lim\limits_{n\rightarrow\omega}k(n)<\infty$. Otherwise, if $\lim\limits_{n\rightarrow\omega}k(n)=\infty$, then since $u_n\in \mathcal U(R_{k(n)})$, for all $n\geq 1$, we would get that $u=(u_n)\in R'\cap R^{\omega}$. Since $\|x-uxu^*\|_2=\lim\limits_{n\rightarrow\omega}\|x_n-u_nx_nu_n^*\|_2\geq\varepsilon/2$, this would contradict that $x$ commutes with $R'\cap  R^{\omega}$.
Now, if  $k=\lim\limits_{n\rightarrow\omega}k(n)+1<\infty$, then $\|x-E_{(R_k'\cap M)^{\omega}}(x)\|_2=\lim\limits_{n\rightarrow\omega}\|x_n-E_{R_{k(n)+1}'\cap M}(x_n)\|_2\leq\varepsilon.$ \hfill$\square$

Finally, note that if $k\geq 1$, then $R_k'\cap M$ is generated by $R'\cap M$ and $R^k$. Since $R^k$ is finite dimensional, it follows that $(R_k'\cap M)^{\omega}$ is generated by $(R'\cap M)^{\omega}$ and $R^k$. In particular, $(R_k'\cap M)^{\omega}\subset (R'\cap M)^{\omega}\vee R$, for every $k\geq 1$. In combination with the claim, this implies that $x\in (R'\cap M)^{\omega}\vee R$, as desired.
\hfill$\blacksquare$

\begin{corollary}\label{commutant}
Let $M=P\bar{\otimes}R$, where $P$ is a tracial von Neumann algebra and $R$ is the hyperfinite II$_1$ factor. 
Then $(R'\cap R^{\omega})'\cap M^{\omega}=P^{\omega}\bar{\otimes}R$.
Moreover, $\mathcal Z(M'\cap M^{\omega})=\mathcal Z(P'\cap P^{\omega})$.
\end{corollary}
{\it Proof.} The first assertion is a consequence of Theorem \ref{comm}. This implies that $$\mathcal Z(M'\cap M^{\omega})\subset (M'\cap M^{\omega})'\cap M^{\omega}\subset (R'\cap R^{\omega})'\cap M^{\omega}= P^{\omega}\bar{\otimes}R,$$ and thus $\mathcal Z(M'\cap M^{\omega})\subset M'\cap (P^{\omega}\bar{\otimes}R)=P'\cap P^{\omega}$. Since clearly $P'\cap P^{\omega}\subset M'\cap M^{\omega}$, we deduce that $\mathcal Z(M'\cap M^{\omega})\subset\mathcal Z(P'\cap P^{\omega})$. Since  $\mathcal Z(P'\cap P^{\omega})\subset \mathcal Z(M'\cap M^{\omega})$ by \cite[Lemma 5.3]{Ma17}, the moreover assertion follows.
\hfill$\blacksquare$

\subsection{A result on normalizers} The next lemma identifies the normalizer of a Cartan subalgebra of a II$_1$ factor $P$ in the tensor product of $P$ with another II$_1$ factor. More generally, we have

\begin{lemma}\label{normalizer}
Let $P\subset M$  be two II$_1$ factors, and $A\subset P$ be a Cartan subalgebra. Assume that $P\vee (A'\cap M)=M$.
If $u\in\mathcal N_{M}(A)$, then we can find $u_1\in\mathcal N_{P}(A)$ and $u_2\in\mathcal U(A'\cap M)$ such that $u=u_1u_2$.
In particular, $\mathcal R(A\subset M)=\mathcal R(A\subset P)$.
\end{lemma}

{\it Proof.}  Let $u\in\mathcal N_{M}(A)$ and $\theta$ be the automorphism of $A$ given by $\theta(a)=uau^*$, for every $a\in A$. 
Let $p\in A$ be a non-zero projection. 

Then $\theta(a)up=upa$, for every $a\in Ap.$ Since $A\subset P$ is a Cartan subalgebra and $P\vee (A'\cap M)=M$, 
the linear span of $\mathcal V=\{wv\mid v\in\mathcal N_{P}(A),w\in\mathcal U(A'\cap M)\}$ is $\|.\|_2$-dense in $M$. 
Thus, we can find $v\in\mathcal N_{P}(A)$ and $w\in\mathcal U(A'\cap M)$ such that $\uptau(upvw)\not=0$, and hence we have $\xi:=E_{A}(upvw)\not=0$. 
Since $\theta(a)upvw=upvw(v^*av)$, by applying $E_{A}$ we conclude that $\theta(a)\xi=(v^*av)\xi$, for all $a\in Ap$.
Thus, if $q\in A$ denotes the support projection of $\xi$, then $\theta(a)q=(v^*av)q$, for all $a\in Ap$. Since $\xi\in A(v^*pv)$, we deduce that $q\in A(v^*pv)$. Thus,  $p'=vqv^*\in Ap$ and we have $\theta(p')q=(v^*p'v)q=q$. Thus, for all $a\in Ap'$ we have that $\theta(a)=\theta(ap')=\theta(a)q=(v^*av)q=v^*av$.

In conclusion, for every non-zero projection $p\in A$ there is a non-zero projection $p'\in Ap$ such that the restriction of $\theta$ to $Ap'$ is equal to $\text{Ad}(v)$, for some $v\in \mathcal N_P(A)$. Since $A\subset P$ is a Cartan subalgebra, this implies that there exists $u_1\in\mathcal N_P(A)$ such that $\theta(a)=u_1au_1^*$, for all $a\in A$. Then $u_1^*u\in A'\cap M$, which proves the conclusion.
\hfill$\blacksquare$

\section{Proofs of Theorem \ref{A} and Corollaries \ref{B} and \ref{C}}
This section is devoted to the proof of our main technical result, Theorem \ref{A}, whose statement we recall below for convenience. We will end the section by deriving Corollaries \ref{B} and \ref{C} from it.
\begin{theorem}\label{Cartantransfer}
Let $P_1$ be a II$_1$ factor which admits a Cartan subalgebra $A_1$ such that $P_1'\cap P_1^{\omega}\subset A_1^{\omega}$. Let $P_2$ be a non-McDuff II$_1$ factor, and    
 $\theta:P_1\bar{\otimes}R_1\rightarrow P_2\bar{\otimes}R_2$ be an isomorphism, where $R_1, R_2$ are hyperfinite II$_1$ factors.

Then $P_2$ admits a Cartan subalgebra $A_2$ satisfying $P_2'\cap P_2^{\omega}\subset A_2^{\omega}$ and we can find a unitary element $u\in P_2\bar{\otimes}R_2$ and $t>0$ such that  $\theta(A_1^t)=uA_2u^*$, where we identify $P_1\bar{\otimes}R_1=P_1^{t}\bar{\otimes}R_1^{1/t}$. 

Moreover, $\mathcal R(A_1\subset P_1)^t$ is isomorphic to $\mathcal R(A_2\subset P_2)$.

\end{theorem}

{\it Proof.} 
Note first that by Lemma \ref{normalizer}, the moreover assertion is a consequence of the main assertion. To prove the main assertion of the theorem,
put $M=P_1\bar{\otimes}R_1$ and identify $M=P_2\bar{\otimes}R_2$ using $\theta$. 
Since $P_1'\cap P_1^{\omega}\subset A_1^{\omega}$, by applying\cite[Proposition 5.2]{Ma17} we deduce that \begin{equation}\label{central1}M'\cap M^{\omega}\subset (A_1\bar{\otimes} R_1)^{\omega}.\end{equation}
Thus, we have that $R_2'\cap R_2^{\omega}\subset (A_1\bar{\otimes}R_1)^{\omega}$. By applying Lemma \ref{intert}, we get that $R_2\prec_{M}A_1\bar{\otimes}R_1$.
Using this fact, the first part of the proof of \cite[Proposition 6.3]{Ho15} provides some $s>0$ such that if we identify $P_2\bar{\otimes}R_2=P_2^{s}\bar{\otimes}R_2^{1/s}$, then we can find a non-zero projection  $q_0\in R_2^{1/s}$ and a unitary $v\in M$ such that $q:=vq_0v^*\in A_1\bar{\otimes}R_1$ and  we have \begin{equation}\label{dan}\text{Ad}(v) (q_0R_2^{1/s}q_0)\subset q(A_1\bar{\otimes}R_1)q.\end{equation}

Denote $P_3=\text{Ad}(v) (P_2^s)$ and $R_3=\text{Ad}(v) (R_2^{1/s})$. Then  $M=P_3\bar{\otimes}R_3$, $q\in R_3$, and \eqref{dan} rewrites as $qR_3q\subset q(A_1\bar{\otimes}R_1)q$. 
By taking relative commutants inside $M$, this gives that $A_1q\subset P_3\otimes q.$
Therefore, we can write $A_1q=A_3\otimes q$, where $A_3\subset P_3$ is a unital von Neumann subalgebra. 

Since $P_2$ is not McDuff, $P_3\cong P_2^s$ is also not McDuff, i.e., $P_3'\cap P_3^{\omega}$ is abelian.  Since $P_1'\cap P_1^{\omega}$ is abelian, Corollary \ref{commutant} implies that $P_3'\cap P_3^{\omega}=\mathcal Z(M'\cap M^{\omega})=P_1'\cap P_1^{\omega}\subset A_1^{\omega}.$
 From this we get that $(P_3'\cap P_3^{\omega})\otimes q\subset A_1^{\omega}q=A_3^{\omega}\otimes q,$ and since $q$ is non-zero, it follows that
\begin{equation}\label{P3} P_3'\cap P_3^{\omega}\subset A_3^{\omega}. \end{equation}
By combining \eqref{P3} with \cite[Proposition 5.2]{Ma17} we derive that \begin{equation}\label{central2}M'\cap M^{\omega}\subset (A_3\bar{\otimes}R_3)^{\omega}.\end{equation}
The rest of the proof is divided between three claims. 

{\bf Claim 1.} $A_3'\cap P_3\prec_{A_3'\cap P_3}A_3$.

{\it Proof of Claim 1.} Since $A_3\bar{\otimes}qR_3q\subset qMq$ is generated by $A_3\otimes  q=A_1q$ and $qR_3q\subset q(A_1\bar{\otimes}R_1)q$, we get that $A_3\bar{\otimes} qR_3q\subset q(A_1\bar{\otimes}R_1)q$. Let us first show that
 $q(A_1\bar{\otimes}R_1)q\prec_{q(A_1\bar{\otimes}R_1)q}A_3\bar{\otimes}qR_3q$. This is equivalent to proving that $(A_1\bar{\otimes}R_1)\tilde q\prec_{(A_1\bar{\otimes}R_1)\tilde q}A_3\bar{\otimes}qR_3q$, where $\tilde q\in\mathcal Z(A_1\bar{\otimes}R_1)=A_1$ denotes the central support of $q$. 
Let $v_1,..,v_m\in A_1\bar{\otimes}R_1$ be partial isometries such that $v_i^*v_i\leq q$, for all $1\leq i\leq m$, the projections $v_iv_i^*, 1\leq i\leq m$, are mutually orthogonal, and $\|\tilde q-\sum_{i=1}^mv_iv_i^*\|_2\leq \|\tilde q\|_2/4.$

Write $R_1=\bar{\otimes}_{k\in\mathbb N}\mathbb M_2(\mathbb C)$, and for $n\in\mathbb N$, define $R_{1,n}=\bar{\otimes}_{k\geq n}\mathbb M_2(\mathbb C)$. Since $\prod_{\omega}R_{1,n}\subset M'\cap M^{\omega}$,
by using \eqref{central2} and Lemma \ref{easy}, we can find $n\in\mathbb N$ such that $\|u-E_{A_3\bar{\otimes}R_3}(u)\|_2\leq\|\tilde q\|_2/(4m)$ and $\|v_iu-uv_i\|_2\leq \|\tilde q\|_2/(4m)$, for all $u\in\mathcal U(R_{1,n})$ and every $1\leq i\leq m$.

Let $a\in\mathcal U(A_1)$ and $u\in\mathcal U(R_{1,n})$. Then we have \begin{align*}\|u\tilde q-\sum_{i=1}^mv_iE_{A_3\bar{\otimes}R_3}(u)v_i^*\|_2&\leq \|u\tilde q-\sum_{i=1}^mv_iuv_i^*\|_2+\|\tilde q\|_2/4\\&\leq \|u\tilde q-u\sum_{i=1}^m v_iv_i^*\|_2+\|\tilde q\|_2/2\\&\leq (3\|\tilde q\|_2)/4.
\end{align*}

Since $a\in A_1=\mathcal Z(A_1\bar{\otimes}R_1)$, we get that $\|au\tilde q-\sum_{i=1}^mv_iaE_{A_3\bar{\otimes}R_3}(u)v_i^*\|_2\leq (3\|\tilde q\|_2)/4<\|\tilde q\|_2.$ Since $v_i\in (A_1\bar{\otimes}R_1)q$, for all $1\leq i\leq m$ and $qa\in A_1q=A_3\otimes q$, we deduce that $$\sum_{i=1}^mv_iaE_{A_3\bar{\otimes}R_3}(u)v_i^*=\sum_{i=1}^mv_i(aq)(qE_{A_3\bar{\otimes}R_3}(u)q)v_i^*\in\sum_{i=1}^mv_i(A_3\bar{\otimes}qR_3q)_1v_i^*.$$
Since the subgroup $\{au\tilde q\mid a\in\mathcal U(A_1),u\in\mathcal U(R_{1,n})\}$ of $\mathcal U((A_1\bar{\otimes}R_{1,n})\tilde q)$ generates $(A_1\bar{\otimes}R_{1,n})\tilde q$, by Theorem \ref{corner} we get that $(A_1\bar{\otimes}R_{1,n})\tilde q\prec_{(A_1\bar{\otimes}R_1)\tilde q}A_3\bar{\otimes}qR_3q$. Since $R_{1,n}\subset R_1$ is a finite index subfactor, we get that $(A_1\bar{\otimes}R_{1})\tilde q\prec_{(A_1\bar{\otimes}R_1)\tilde q}A_3\bar{\otimes}qR_3q$, and thus $q(A_1\bar{\otimes}R_1)q\prec_{q(A_1\bar{\otimes}R_1)q}A_3\bar{\otimes}qR_3q$. 
Since $q(A_1\bar{\otimes}R_1)q=(A_1q)'\cap qMq=(A_3\otimes q)'\cap qMq=(A_3'\cap P_3)\bar{\otimes}qR_3q$, the last embedding rewrites as $(A_3'\cap P_3)\bar{\otimes}qR_3q\prec_{(A_3'\cap P_3)\bar{\otimes}qR_3q}A_3\bar{\otimes}qR_3q$, which readily implies the claim.
\hfill$\square$

{\bf Claim 2.} There exist a Cartan subalgebra $A_4\subset P_3$ such that $P_3'\cap P_3^{\omega}\subset A_4^{\omega}$, and a non-zero projection $p\in A_4$ such that $r=p\otimes q$ belongs to $q(A_1\bar{\otimes}R_1)q$ and $A_4r=A_1r$.

{\it Proof of Claim 2.}
By Claim 1 we can find a non-zero projection $p\in A_3'\cap P_3$ such that $A_3p$ is maximal abelian in $pP_3p$.
Let $r=p\otimes q$. Then $r\in (A_3\otimes q)'\cap qMq=(A_1q)'\cap qMq=q(A_1\bar{\otimes}R_1)q$ and $A_3r=(A_3\otimes q)r=(A_1q)r=A_1r$. Since $A_1\subset P_1$ is a Cartan subalgebra, $A_1$ is regular in $M=P_1\bar{\otimes}R_1$. By \cite[Lemma 3.5.1]{Po03} we get that $A_3r=A_1r$ is quasi-regular in $rMr$. This implies that $A_3r=A_3p\otimes q$ is quasi-regular in 
$r(P_3\otimes q)r=pP_3p\otimes q$, hence $A_3p$ is quasi-regular in $pP_3p$.
As $A_3p\subset pP_3p$ is maximal abelian, \cite[Theorem 2.7]{PS03} gives that $A_3p\subset pP_3p$ is a Cartan subalgebra.

Since  $P_3$ is a II$_1$ factor, we can find a Cartan subalgebra $A_4\subset P_3$ such that $p\in A_4$ and $A_4p=A_3p$. 
Thus, we have that $A_4r=A_3r=A_1r$. In order to finish the proof of the claim, it remains to argue that  $P_3'\cap P_3^{\omega}\subset A_4^{\omega}$. 
First, by \eqref{P3}, we get that $(pP_3p)'\cap (pP_3p)^{\omega}\subset (A_3p)^{\omega}=(A_4p)^{\omega}$.
Since $P_3$ is a II$_1$ factor and $A_4$ is a Cartan subalgebra, we can find partial isometries $v_1,...,v_k\in P_3$ such that  $v_i^*v_i\in A_4p$ and $v_iA_4v_i^*\subset A_4$, for every $1\leq i\leq k$, and $\sum_{i=1}^kv_iv_i^*=1$.  If $x\in P_3'\cap P_3^{\omega}$, then 
$v_i^*xv_i\in (pP_3p)'\cap (pP_3p)^{\omega}\subset (A_4p)^{\omega}$ and thus $v_iv_i^*xv_iv_i^*=v_i(v_i^*xv_i)v_i^*\in A_4^{\omega}$, for every $1\leq i\leq k$. Since $x=\sum_{i=1}^kv_iv_i^*xv_iv_i^*$, we get that $x\in A_4^{\omega}$, as claimed. 
\hfill$\square$

Since $A_4p$ is a Cartan subalgebra of $pP_3p$, we deduce that $A_1r=A_4r=A_4p\otimes q$ is a regular subalgebra of $rMr$. 
 Using this fact, we will next show the following:

{\bf Claim 3.} There exist a projection $e\in A_1$ and $\alpha\geq 0$ such that $E_{A_1}(r)=\alpha e$.

{\it Proof of Claim 3.} 
Denote by $e\in A_1$ the support projection of $E_{A_1}(r)$.  Let $u\in\mathcal N_{rMr}(A_1r)$.  Since
 $A_1e\ni a\mapsto ar\in A_1r$ is an isomorphism, we can find an isomorphism $\theta:A_1e\rightarrow A_1e$ such that \begin{equation}\label{theta}uaru^*=\theta(a)r,\;\;\;\text{for all $a\in A_1e$.}\end{equation}
 Thus, $\theta(a)=E_{A_1}(uaru^*)E_{A_1}(r)^{-1}$, for all $a\in A_1e$. Combining this with the fact that $A_1e\subset eP_1e$ is a Cartan subalgebra, we can find $v\in\mathcal N_{eP_1e}(A_1e)$ such that $\theta(a)=vav^*$, for all $a\in A_1e$. Thus, $\uptau(av^*rv)=\uptau(vav^*r)=\uptau(\theta(a)r)=\uptau(uaru^*)=\uptau(ar)$, for all $a\in A_1e$, hence $E_{A_1}(v^*rv)=E_{A_1}(r)$. Since $v$ normalizes $A_1e$, we derive that $E_{A_1}(r)$ commutes with $v$. Since $E_{A_1}(r)\in A_1e$ we get that $\theta(E_{A_1}(r))=E_{A_1}(r)$. In combination with \eqref{theta}, we get that $E_{A_1}(r)r=\theta(E_{A_1}(r))r=uE_{A_1}(r)ru^*.$
This shows that $E_{A_1}(r)r\in rMr$ commutes with every $u\in\mathcal N_{rMr}(A_1r)$. Since $A_1r$ is regular in $rMr$ and $rMr$ is a  factor, we derive that $E_{A_1}(r)r=\alpha r$, for some $\alpha\geq 0$. Hence $E_{A_1}(r)^2=\alpha E_{A_1}(r)$, which implies the claim.\hfill$\square$

Let $f\in R_1$ be a projection such that $\uptau(f)=\alpha$. By Claim  3 we have $E_{A_1}(r)=E_{A_1}(e\otimes f)$, thus we can find a unitary $w\in A_1\bar{\otimes}R_1$ such that $r=w(e\otimes f)w^*$. Altogether, we have found a Cartan subalgebra $A_4\subset P_3$, projections $p\in A_4, q\in R_3,e\in A_1, f\in R_1$, and a unitary $w\in M$ such that \begin{equation}\label{conjugacy} A_4p\otimes q=A_4r=A_1r=\text{Ad}(w)(A_1e\otimes f).\end{equation}
If $P$ and $R$ are II$_1$ factors, $A\subset P$ is a Cartan subalgebra, and $a\in A, b\in R$ are projections, then we can identify $P\bar{\otimes}R=P^{1/\uptau(b)}\bar{\otimes}R^{\uptau(b)}$ such that $Aa\otimes b$ is identified with $A^{1/\uptau(b)}c$, where $c\in A^{1/\uptau(b)}\subset P^{1/\uptau(b)}$ is a projection of trace  $\uptau(a)\uptau(b)$. 
By combining this fact with \eqref{conjugacy} it follows that we can identify $M=P_1^t\bar{\otimes}R_1^{1/t}$, for $t=\uptau(q)/\uptau(f)$, such that $A_4$ and $A_1^t$ are unitarily conjugate. 

Finally, recall that $P_3=\text{Ad}(v) (P_2^s)$ and $R_3=\text{Ad}(v) (R_2^{1/s})$.  Let  $A_2\subset P_2$ be a Cartan subalgebra such that $A_4=\text{Ad}(v)(A_2^s)$. 
 Since $P_3'\cap P_3^{\omega}\subset A_4^{\omega}$ by Claim 2, we get that $(P_2^s)'\cap (P_2^s)^{\omega}\subset (A_2^s)^{\omega}$ and the argument from the proof of Claim 2 implies that $P_2'\cap P_2^{\omega}\subset A_2^{\omega}$. Moreover, identifying  $M=P_1^t\bar{\otimes}R_1^{1/t}=P_2^s\bar{\otimes}R_2^{1/s}$ we have that $A_1^t$ and $A_2^s$ are unitarily conjugate. Thus, if we identify $M=P_1^{t/s}\bar{\otimes}R_1^{s/t}$, then $A_1^{t/s}$ is unitarily conjugate to $A_2$. This finishes the proof.
\hfill$\blacksquare$

\subsection*{Proof of Corollary \ref{B}} Let $M=L^{\infty}(X)\rtimes\mathbb F_n$, where $\mathbb F_n\curvearrowright (X,\mu)$ is a free ergodic p.m.p. action, for some $n\geq 2$. Let $N$ be a II$_1$ factor such that $M\bar{\otimes}R\cong N\bar{\otimes}R$.
If $N$ is McDuff, then $N\cong N\bar{\otimes}R\cong M\bar{\otimes}R$. Thus, we may assume that $N$ is not McDuff.

 Since $\mathbb F_n$ is not inner amenable, \cite{Ch82} implies that $M'\cap M^{\omega}\subset L^{\infty}(X)^{\omega}$. By Theorem  \ref{Cartantransfer} we can find a Cartan subalgebra $B\subset N$, some $t>0$, and a unitary $u\in N\bar{\otimes}R$ such that $\theta(L^{\infty}(X)^t)=uBu^*$, where we identify $M\bar{\otimes}R=M^t\bar{\otimes}R^{1/t}$. 
Moreover, if we put  $\mathcal R=\mathcal R(B\subset N)$, then $\mathcal R\cong\mathcal R(\mathbb F_n\curvearrowright X)^t.$
Since $\mathcal R(\mathbb F_n\curvearrowright X)^t$ is a treeable equivalence relation,  $\mathcal R$ is also treeable. Since $\mathcal R$ is treeable, we have that H$^2(\mathcal R,\mathbb T)=0$ (see e.g. \cite[Corollary 2.4]{Ki14}). Therefore, \cite{FM77} implies that $N\cong L(\mathcal R)$ and thus $N\cong L(\mathcal R(\mathbb F_n\curvearrowright X)^t)\cong M^t$, as desired. \hfill$\blacksquare$

\subsection*{Proof of Corollary \ref{C}} For $i\in\{1,2\}$, put $A_i=L^{\infty}(X_i)$ and $M_i=L(\mathcal R_i)$. Denote $B=L^{\infty}(Y)$ and $R=L(\mathcal T)$. Since $\mathcal T$ is hyperfinite, ergodic and p.m.p., $R$ is a hyperfinite II$_1$ factor.  Since $\mathcal R_1\times\mathcal T\cong \mathcal R_2\times\mathcal T$, we have an isomorphism $\theta:M_1\bar{\otimes}R\rightarrow M_2\bar{\otimes}R$ such that $\theta(A_1\bar{\otimes} B)=A_2\bar{\otimes}B$. Since $M_1'\cap M_1^{\omega}\subset A_1^{\omega}$ and $M_2$ is not McDuff, by applying Theorem \ref{A} we deduce the existence of a unitary $u\in M_2\bar{\otimes}R$, a Cartan subagebra $\tilde A_2$ of $M_2$, and some $t>0$,  such that $\theta(A_1^t)=u\tilde A_2u^*$, where we identify the inclusions $(A_1\bar{\otimes}B\subset M_1\bar{\otimes}R)$ and $(A_1^t\bar{\otimes}B^{1/t}\subset M_1^t\bar{\otimes}R^{1/t})$.

We claim that $\tilde A_2\prec_{M_2}A_2$. Assume by contradiction that $\tilde A_2\not\prec_{M_2}A_2$. Then we can find a sequence $v_n\in \mathcal U(\tilde A_2)$ such that $\|E_{A_2}(xv_ny)\|_2\rightarrow 0$, for all $x,y\in M_2$. Then $\|E_{A_2\bar{\otimes}B}(zv_nt)\|_2\rightarrow 0$, for all $z,t\in M_2\bar{\otimes}R$. Indeed, it is enough to check this in the case when $z=z_1\otimes z_2$ and $t=t_1\otimes t_2$, for some $z_1,t_1\in M_2, z_2,t_2\in R$. In this case, we have that $E_{A_2\bar{\otimes}B}(zv_nt)=E_{A_2}(z_1v_nt_1)\otimes z_2t_2$, and hence $\|E_{A_2\bar{\otimes}B}(zv_nt)\|_2=\|E_{A_2}(z_1v_nt_1)\|_2\|z_2t_2\|_2\rightarrow 0$. Since $\theta(A_1^{t}\bar{\otimes}B^{1/t})=A_2\bar{\otimes}B$, we get that $\|E_{A_1^t\bar{\otimes}B^{1/t}}(\theta^{-1}(uv_nu^*))\|_2=\|E_{A_2\bar{\otimes}B}(uv_nu^*)\|_2\rightarrow 0.$ This however contradicts the fact that $\theta^{-1}(uv_nu^*)\in A_1^t\subset A_1^t\bar{\otimes}B^{1/t}$, thus proving the claim.

Since $\tilde A_2$ and $A_2$ are Cartan subalgebras of $M_2$, by using the above claim and \cite[Theorem A.1]{Po01}, we get that $\tilde A_2=vA_2v^*$, for some $v\in\mathcal U(M_2)$. Thus,  $\theta(A_1^t)=\text{Ad}(u)(\tilde A_2)=\text{Ad}(uv)(A_2)$. Finally, Lemma \ref{normalizer} gives that 
$$\mathcal R_1^t\cong\mathcal R(A_1^t\subset M_1^t\bar{\otimes}R^{1/t})\cong\mathcal R(\theta(A_1^t)\subset M_2\bar{\otimes}R)\cong\mathcal R(A_2\subset M_2\bar{\otimes}R)\cong\mathcal R_2,$$ which proves the conclusion.
\hfill$\blacksquare$

\section{A cocycle rigidity result}
This section is devoted to the proof of the following cocycle rigidity result. This result, which seems to be of independent interest, will be used later on to derive Theorems \ref{D} and \ref{H}. 
\begin{theorem}\label{rigid}
Let $\mathcal S$ and $\mathcal T$ be countable ergodic p.m.p. equivalence relations on probability spaces $(X,\mu)$ and $(Y,\nu)$. Let $\pi:(X,\mu)\rightarrow (Y,\nu)$ be a factor map, consider the embedding $L^{\infty}(Y)\subset L^{\infty}(X)$ given by $f\mapsto f\circ\pi$, and denote $\mathcal S_0=\{(x_1,x_2)\in \mathcal S\mid \pi(x_1)=\pi(x_2)\}$.  Assume that \begin{enumerate}
\item\label{1} $[\pi(x)]_{\mathcal T}\subset\pi([x]_{\mathcal S})$, for almost every $x\in X$,
\item\label{2} $L^{\infty}(X)^{[\mathcal S_0]}\subset L^{\infty}(Y)$, and 
\item\label{3} $(L^{\infty}(X)^{\omega})^{[\mathcal S]}\subset L^{\infty}(Y)^{\omega}$.
\end{enumerate}

Let $\alpha:X\rightarrow\text{Aut}(R)$ be a measurable map satisfying $\alpha(x)^{-1}\alpha(y)\in \text{Inn}(R)$, for almost every $(x,y)\in\mathcal S$, where $R$ denotes the hyperfinite II$_1$ factor.

Then we can find measurable maps  $\tilde\alpha:Y\rightarrow\text{Aut}(R)$ and $U:X\rightarrow\mathcal U(R)$ such that we have $\alpha(x)=\tilde\alpha(\pi(x))\text{Ad}(U(x))$, for almost every $x\in X$. 

\end{theorem}

{\it Proof.} In order to prove the conclusion, we may assume that $\alpha(x)^{-1}\alpha(y)\in\text{Inn}(R)$, for all $(x,y)\in\mathcal S$. Consider the disintegration $\mu=\int_{Y}\mu_{y}\text{d}\nu(y)$, where $\mu_{y}$ is a probability measure on $X$ supported on $\pi^{-1}(\{y\})$, for all $y\in Y$. Put $\tilde X=\{(x_1,x_2)\in X\times X\mid \pi(x_1)=\pi(x_2)\}$ and endow $\tilde X$ with the probability measure $\tilde\mu=\int_{Y}(\mu_y\times\mu_y)\;\text{d}\nu(y).$ Let $A=\{(x_1,x_2)\in\tilde X\mid \alpha(x_1)^{-1}\alpha(x_2)\in\text{Inn}(R)\}$. 

The first part of the proof consists of using condition \eqref{3} to establish the following:

{\bf Claim.} $\tilde\mu(A)>0$.

{\it Proof of the claim.} Denote by $L(\mathcal S)\subset \mathbb B(L^2(\mathcal S))$ the von Neumann algebra associated to $\mathcal S$ and by $\{u_{\gamma}\}_{\gamma\in [\mathcal S]}\subset L(\mathcal S)$ the canonical unitaries \cite{FM77}.
Condition \eqref{3} gives that $L^{\infty}(X)^{\omega}\cap L(\mathcal S)'\subset L^{\infty}(Y)^{\omega}$. By \cite[Proposition 5.2]{Ma17} we deduce that \begin{equation}\label{incl1} (L^{\infty}(X)\bar{\otimes}R)^{\omega}\cap (L(\mathcal S)\otimes 1)'\subset (L^{\infty}(Y)\bar{\otimes}R)^{\omega}.\end{equation}

Next, we identify $L^{\infty}(X)\bar{\otimes}R=L^{\infty}(X,R)$, and define  $\theta$ to be the automorphism of  $L^{\infty}(X)\bar{\otimes}R$ given by $\theta(T)(x)=\alpha(x)(T(x))$, for every bounded measurable function $T:X\rightarrow R$ and all $x\in X$.

Fix $\gamma\in [\mathcal S]$, and define $w(\gamma,x)=\alpha(x)^{-1}\alpha(\gamma^{-1}x)\in\text{Inn}(R)$, for $x\in X$. If $T\in R$, then we have $\text{Ad}(u_{\gamma}\otimes 1)(\theta(1\otimes T))(x)=\alpha(\gamma^{-1}x)(T)=\alpha(x)(w(\gamma,x)(T))$. 
Thus, if $(T_n)_n\in R'\cap R^{\omega}$, then as $w(\gamma,x)\in\text{Inn}(R)$ we get that $\|\text{Ad}(u_{\gamma}\otimes 1)(\theta(1\otimes T_n))(x)-\theta(1\otimes T_n)(x)\|_2=\|w(\gamma,x)(T_n)-T_n\|_2\rightarrow 0$, for almost every $x\in X$. Thus,  $\|\text{Ad}(u_{\gamma}\otimes 1)(\theta(1\otimes T_n))-\theta(1\otimes T_n)\|_2\rightarrow 0$, by the dominated convergence theorem. 
As this holds for all $\gamma\in [\mathcal S]$, we get that 
$$\theta(1\otimes R)'\cap \theta(1\otimes R)^{\omega}
\subset (L(\mathcal S)\otimes 1)'\cap (L^{\infty}(X)\bar{\otimes}R)^{\omega}.$$
By combining this with \eqref{incl1}, we get that $\theta(1\otimes R)'\cap \theta(1\otimes R)^{\omega}\subset (L^{\infty}(Y)\bar{\otimes}R)^{\omega}$. Lemma \ref{intert} then gives that $\theta(1\otimes R)\prec_{L^{\infty}(X)\bar{\otimes}R}L^{\infty}(Y)\bar{\otimes}R.$
Thus, we can find a  projection $p\in R$, a $*$-homomorphism $\rho:pRp\rightarrow L^{\infty}(Y)\bar{\otimes}R$, and a non-zero partial isometry $v\in \rho(p)(L^{\infty}(X)\bar{\otimes} R)\theta(1\otimes p)$ such that \begin{equation}\label{rho1}\rho(T)v=v\theta(1\otimes T),\;\;\;\text{for all}\;\;\;T\in pRp.\end{equation} 
If we view $\rho(T)$ as a function $\rho(T):Y\rightarrow R$, then \eqref{rho1} rewrites as 
\begin{equation} \label{rho2}
\rho(T)(\pi(x))v(x)=v(x)\alpha(x)(T),\;\;\;\text{for all $T\in pRp$ and almost every $x\in X$.}\end{equation}

Since $v(x)v(x)^*\leq\rho(p)(\pi(x))$, for almost every $x\in X$, \eqref{rho2} gives that for all $T\in\mathcal U(pRp)$, we have
\begin{equation}\label{rho3}\alpha(x_1)(T)v(x_1)^*v(x_2)=v(x_1)^*v(x_2)\alpha(x_2)(T),\;\;\;\text{for  almost every $(x_1,x_2)\in \tilde X.$}
\end{equation}

Now, note first that if $\alpha_1,\alpha_2\in\text{Aut}(R)$ satisfy $\alpha_1(T)w=w\alpha_2(T)$ for all $T\in pRp$, and for some non-zero partial isometry $w\in \alpha_1(p)R\alpha_2(p)$, then $\alpha_1^{-1}\alpha_2\in\text{Inn}(R)$. 
Second, note that \begin{align*}
\int_{\tilde X}\|v(x_1)^*v(x_2)\|_2^2\;\text{d}\tilde\mu(x_1,x_2)&=\int_{Y}\Big(\int_{X\times X}\uptau(v(x_1)v(x_1)^*v(x_2)v(x_2)^*)\;\text{d}(\mu_y\times\mu_y)(x_1,x_2)\Big)\;\text{d}\nu(y)\\&=\int_{Y}\uptau(E_{L^{\infty}(Y)\bar{\otimes}R}(vv^*)(y)^2)\;\text{d}\nu(y)\\&=\|E_{L^{\infty}(Y)\bar{\otimes}R}(vv^*)\|_2^2>0.
\end{align*}
By combining the last two facts and \eqref{rho3}, the claim follows. \hfill$\square$

In the second part of the proof, we  will use the ergodicity assumptions on $\mathcal S_0$ and $\mathcal T$ together with the claim to deduce the conclusion. First, we have that $\mu_y$ is $\mathcal S_0$-invariant, while condition \eqref{2} implies that $\mu_y$ is $\mathcal S_0$-ergodic, for almost every $y\in Y$. If $\gamma\in [\mathcal S_0]$ and $(x_1,x_2)\in A$, then since $\pi(\gamma x_1)=\pi(x_1)$ and $\alpha(\gamma x_1)^{-1}\alpha(x_1)\in\text{Inn}(R)$, we derive that $(\gamma x_1,x_2)\in A$. Thus, for all $\gamma\in [\mathcal S_0]$, the set $A^{x_2}=\{x_1\in\pi^{-1}(\{\pi(x_2)\})\mid (x_1,x_2)\in A\}$ satisfies $\mu_{\pi(x_2)}(\gamma A^{x_2}\Delta A^{x_2})=0$, for almost every $x_2\in X$. Since $\mu_{\pi(x_2)}$ is $\mathcal S_0$-ergodic, we get that $\mu_{\pi(x_2)}(A^{x_2})\in\{0,1\}$, for almost every $x_2\in X$. Similarly, we get that $\mu_{\pi(x_1)}(A_{x_1})\in\{0,1\}$, for almost every $x_1\in X$, where $A_{x_1}=\{x_2\in\pi^{-1}(\{\pi(x_1)\})\mid (x_1,x_2)\in A\}$. 
Since $\tilde\mu(A)>0$, the last two facts imply the existence of a measurable set $B\subset Y$ such that $\nu(B)>0$ and $A=\{(x_1,x_2)\in \tilde X\mid \pi(x_1)=\pi(x_2)\in B\}$.

We claim that $B$ is $\mathcal T$-invariant. Let $\gamma\in [\mathcal T]$. Since $\gamma y\in [y]_{\mathcal T}$, by condition \eqref{1}, for almost every $y\in Y$, for  $(\mu_{\gamma y}\times\mu_{\gamma y})$-almost every $(x_1',x_2')\in\pi^{-1}(\{\gamma y\})\times\pi^{-1}(\{\gamma y\})$, we can find $x_1\in [x_1']_{\mathcal S}$ and  $x_2\in [x_2']_{\mathcal S}$ such that $\pi(x_1)=\pi(x_2)=y$.
Since $(x_i,x_i')\in\mathcal S$, $\alpha(x_i)^{-1}\alpha(x_i')\in\text{Inn}(R)$, for $i\in\{1,2\}$.  Additionally, for almost every $y\in B$, we have that $\alpha(x_1)^{-1}\alpha(x_2)\in\text{Inn}(R)$, for $(\mu_y\times\mu_y)$-almost every $(x_1,x_2)\in\pi^{-1}(\{y\})\times\pi^{-1}(\{y\})$. Combining these facts gives that for almost every $y\in B$, we have $\alpha(x_1')^{-1}\alpha(x_2')\in\text{Inn}(R)$, for  $(\mu_{\gamma y}\times\mu_{\gamma y})$-almost every $(x_1',x_2')\in\pi^{-1}(\{\gamma y\})\times\pi^{-1}(\{\gamma y\})$. 
This shows that $\nu(\gamma B\Delta B)=0$. Since $\gamma\in[\mathcal T]$ is arbitrary, we derive that $B$ is indeed $\mathcal T$-invariant. 

Since $\mathcal T$ is ergodic, we conclude that $B=Y$ and thus $A=\tilde X$, almost everywhere. In other words, $\alpha(x_1)^{-1}\alpha(x_2)\in\text{Inn}(R)$, for almost every $(x_1,x_2)\in\tilde X$. 
Let $C=\{y\in Y\mid\mu_y\;\text{is a non-atomic measure}\}$. Denote by $\text{Leb}$ the Lebesgue measure on $[0,1]$.
Then we can find a measure space isomorphism $\sigma: (C\times [0,1],\nu_{|C}\times\text{Leb})\rightarrow (\pi^{-1}(C),\mu_{|\pi^{-1}(C)})$ satisfying $\sigma(\{y\}\times [0,1])=\pi^{-1}(\{y\})$, for all $y\in C$. Since $\alpha(\sigma(y,t_1))^{-1}\alpha(\sigma(y,t_2))\in\text{Inn}(R)$, for $(\nu\times\text{Leb}\times\text{Leb})$-almost every $(y,t_1,t_2)\in C\times [0,1]\times [0,1]$,  Fubini's theorem implies the existence of $t_1\in [0,1]$ such that $\alpha(\sigma(y,t_1))^{-1}\alpha(\sigma(y,t_2))\in\text{Inn}(R)$, for $(\mu\times\text{Leb})$-almost every $(y,t_2)\in C\times [0,1]$. Since $\mu_y$ has atoms for every $y\in Y\setminus{C}$, we can find a measurable map $\eta:Y\setminus {C}\rightarrow X$ such that $\eta(y)\in\pi^{-1}(\{y\})$ and $\mu_y(\{\eta(y)\})>0$, for all $y\in Y\setminus C$. It is now easy to see that if we define $\tilde\alpha:Y\rightarrow\text{Aut}(R)$  by $$\tilde\alpha(y)=\begin{cases}\alpha(\sigma(y,t_1))\;\;\;\text{if $y\in C$, and}\\ \alpha(\eta(y))\;\;\;\;\text{if $y\in Y\setminus C$} \end{cases}$$ then $\beta(x):=\tilde\alpha(\pi(x))^{-1}\alpha(x)\in\text{Inn}(R)$, for almost every $x\in X$.

Finally, let $\mathcal U\subset\mathcal U(R)$ be a countable $\|.\|_2$-dense set. Then an automorphism $\theta$ of $R$ is inner if and only if $\inf_{u\in\mathcal U}\big(\sup_{x\in\mathcal U}\|\theta(x)-uxu^*\|_2\big)=0$. This implies that $\text{Inn}(R)$ is a Borel subset of Aut$(R)$. 
Endow $\text{Inn}(R)$ with the Borel measure $\beta_*\mu$ and consider the Borel map $\zeta:\mathcal U(R)\rightarrow \text{Inn}(R)$ given by $\zeta(u)=\text{Ad}(u)$. 
Since $\zeta$ is onto, by applying \cite[Theorem A.16]{Ta01} we can find a $\beta_*\mu$-measurable map $\xi:\text{Inn}(R)\rightarrow\mathcal U(R)$ such that $\zeta(\xi(\theta))=\theta$, for all $\theta\in\text{Inn}(R)$. After modifying $\xi$ on a $\beta_*\mu$-null set, we can find a Borel map $\xi':\text{Inn}(R)\rightarrow\mathcal U(R)$ such that $\zeta(\xi'(\theta))=\theta$, for $\beta_*\mu$-almost every $\theta\in\text{Inn}(R)$.
Then the map $U:X\rightarrow\mathcal U(R)$ given by $U(x)=\xi'(\beta(x))$ is measurable and satisfies $\tilde\alpha(\pi(x))\text{Ad}(U(x))=\tilde\alpha(\pi(x))\beta(x)=\alpha(x)$, for almost every $x\in X$, which proves the theorem.
\hfill$\blacksquare$

\section{The Jones-Schmidt property}

The main goal of this section is to prove Theorem \ref{D}. We start by analyzing closer the structure of isomorphisms which satisfy the conclusion of Theorem \ref{A}. 

\begin{proposition}\label{cocycle}
Let $M,N,P$ be II$_1$ factors. Assume that $A\subset M$ and  $B\subset N$ are Cartan subalgebras. Denote $\mathcal S=\mathcal R(A\subset M)$, and identify $A=L^{\infty}(X)$ and $B=L^{\infty}(Z)$, for some probability spaces $(X,\mu)$ and $(Z,\eta)$. For $g\in[\mathcal S]$, define $\sigma_g(a)=a\circ g^{-1}$, for every $a\in A$, and let $u_g\in\mathcal N_{M}(A)$ be such that $u_gau_g^*=\sigma_g(a)$, for every $a\in A$.

Let  $\theta:M\bar{\otimes}P\rightarrow N\bar{\otimes}P$ be an isomorphism such that $\theta(A)=B$. Let $\alpha:(X,\mu)\rightarrow (Z,\eta)$ be the measure space isomorphism given by $\theta(a)=a\circ\alpha^{-1}$, for every $a\in A$.

Then $\theta(A\bar{\otimes}P)=B\bar{\otimes}P$, and we can find $w_g\in \mathcal U(A\bar{\otimes}P)$ and $v_g\in\mathcal N_{N}(B)$ for every $g\in [\mathcal S]$, and a measurable map $\varphi:X\rightarrow\text{Aut}(P)$ such that the following conditions hold:

\begin{enumerate}
\item $\theta(w_gu_g)=v_g$ and $\theta(w_g(\sigma_g\otimes\text{id}_P)(w_h)w_{gh}^*)=v_gv_hv_{gh}^*\in B$, for all $g,h\in [\mathcal S]$.
\item $\theta^{-1}(T)(x)=\varphi_{x}(T(\alpha(x)))$, for every $T\in B\bar{\otimes}P$ and almost every $x\in X$, where we identify $A\bar{\otimes}P=L^{\infty}(X,P)$ and $B\bar{\otimes}P=L^{\infty}(Z,P)$.
\item $\text{Ad}(w_g(x))=\varphi_x\varphi_{g^{-1}x}^{-1}$, for all $g\in [\mathcal S]$ and almost every $x\in X$.
\end{enumerate}
\end{proposition}

{\it Proof.} Since $\theta(A)=B$, we get that $\theta(A\bar{\otimes}P)=\theta(A)'\cap (N\bar{\otimes}P)=B'\cap (N\bar{\otimes}P)=B\bar{\otimes}P$.
If $g\in [\mathcal S]$, then $\theta(u_g)\in\mathcal N_{N\bar{\otimes}P}(B)$, hence by Lemma \ref{normalizer} we can find $z_g\in \mathcal U(B\bar{\otimes}P)$ and $v_g\in\mathcal N_{N}(B)$ such that $\theta(u_g)=z_gv_g$. Then $w_g:=\theta^{-1}(z_g^*)\in\mathcal U(A\bar{\otimes}P)$, and we have that $\theta(w_gu_g)=v_g$. Moreover, we get that $v_gv_hv_{gh}^*=\theta(w_gu_gw_hu_hu_h^*u_g^*w_{gh}^*)=\theta(w_g(\sigma_g\otimes\text{id}_{P})(w_h)w_{gh}^*)\in N\cap (B\bar{\otimes}P)=B$, proving (1). 

Since $\theta^{-1}_{|(B\bar{\otimes}P)}:B\bar{\otimes}P\rightarrow A\bar{\otimes}P$ is a trace preserving isomorphism such that $\theta^{-1}(b)=b\circ\alpha$, for every $b\in B$, \cite[Theorem~IV.8.23]{Ta01} provides the existence of a measurable map $\varphi:X\rightarrow\text{Aut}(P)$ satisfying (2). 

Let $T\in P\subset N\bar{\otimes}P$ and $g\in [\mathcal S]$. Since $v_g\in N$, we have that $T$ commutes with $v_g$ and therefore $\theta^{-1}(T)=\theta^{-1}(v_gTv_g^*)=w_gu_g\theta^{-1}(T)u_g^*w_g^*$.  On the other hand, by (2) we get $\theta^{-1}(T)(x)=\varphi_x(T)$, for almost every $x\in X$. Combining these two facts implies that $\varphi_x(T)=w_g(x)\varphi_{g^{-1}x}(T)w_g(x)^*$, for almost every $x\in X$. Since this holds for all $T\in P$ and $g\in [\mathcal S]$,  it follows that (3) holds.
\hfill$\blacksquare$

We continue with the following vanishing cohomology result. 
\begin{lemma}\label{trivial}
Let $\mathcal T$ be a  hyperfinite p.m.p. equivalence relation on a probability space $(X,\mu)$, and $R$ be the hyperfinite II$_1$ factor. Let $c_n:\mathcal T\rightarrow\mathcal U(R)$ be a measurable cocycle, for every $n\geq 1$. Assume that $\|\text{Ad}(c_n(x,y))(v)-v\|_2\rightarrow 0$, as $n\rightarrow\infty$, for every $v\in R$, and almost every $(x,y)\in\mathcal T$.

Then we can find a subsequence $\{n_k\}_{k\geq 1}$ of $\mathbb N$ and  measurable maps $d_k:X\rightarrow\mathcal U(R)$, for all $k\geq 1$, such that $\|\text{Ad}(d_k(x))(v)-v\|_2\rightarrow 0$ and $\|c_{n_k}(x,y)-d_k(x)d_k(y)^{-1}\|_2\rightarrow 0$, as $k\rightarrow\infty$, for every $v\in R$, and almost every $(x,y)\in\mathcal T$.
\end{lemma}

{\it Proof.}
 Let $\{\mathcal T_k\}_{k\geq 1}\subset\mathcal T$ be an increasing sequence of finite subequivalence relations such that we have $\mathcal T=\cup_{k\geq 1}\mathcal T_k$. Moreover, we may assume that $\sup_{x\in X}\#[x]_{\mathcal T_k}<\infty$, for all $k\in\mathbb N$.  For $k\geq 1$, let $\psi_k:X\rightarrow X$ be a measurable map such that $\psi_k(x)\in [x]_{\mathcal T_k}$ and $\psi_k(x)=\psi_k(y)$, for every $(x,y)\in\mathcal T_k$.
Write $R=\bar{\otimes}_{n\geq 1}\mathbb M_2(\mathbb C)$ and define $R_k=\bar{\otimes}_{n\geq k}\mathbb M_2(\mathbb C)$, for all $k\geq 1$. 
For $n,k\geq 1$, define $$X_{n,k}=\{x\in X\mid \|c_n(x,y)-E_{R_k}(c_n(x,y))\|_2\leq 1/k,\;\text{for all $y\in [x]_{\mathcal T_k}$}\}.$$

Then for every $k\geq 1$ we have  $\mu(X_{n,k})\rightarrow 1$, as $n\rightarrow\infty$. Thus,  for all $k\geq 1$, we can find $n_k\geq 1$, such that $\mu(X_{n_k,k})\geq 1-1/2^k$. For $k\geq 1$, we define $d_k:X\rightarrow\mathcal U(R)$ by letting $d_k(x)=c_{n_k}(x,\psi_{k}(x))$. Then for all $k\geq 1$ and $x\in X_{n_k,k}$ we have $\|d_k(x)-E_{R_k}(d_k(x))\|_2\leq 1/k$. The Borel-Cantelli lemma implies that $\|d_k(x)-E_{R_k}(d_k(x))\|_2\rightarrow 0$, for almost every $x\in X$.  Moreover, for all $k\geq 1$ and $(x,y)\in\mathcal T_k$ we have $c_{n_k}(x,y)=d_k(x)d_k(y)^{-1}.$ Since $\mathcal T=\cup_{k\geq 1}\mathcal T_k$, the conclusion follows.
\hfill$\blacksquare$

We next show that the Jones-Schmidt property is equivalent with a formally stronger property.

\begin{proposition}\label{eqJS}  Let $\mathcal S$ be a countable ergodic p.m.p. equivalence relation on a probability space $(X,\mu)$ with the Jones-Schmidt property.
Then there exist a hyperfinite ergodic p.m.p. equivalence relation $\widetilde{\mathcal T}$ on a probability space $(\widetilde Y,\tilde\nu)$ and a factor map 
$\widetilde\pi:(X,\mu)\rightarrow (\widetilde Y,\tilde\nu)$ such that

\begin{enumerate}[label=(\alph*)]
\item $\widetilde\pi([x]_{\mathcal S})=[\widetilde\pi(x)]_{\widetilde{\mathcal T}}$, for almost every $x\in X$, and
\item $(L^{\infty}(X)^{\omega})^{[\mathcal S_0]}=L^{\infty}(\widetilde Y)^{\omega}$, where $\mathcal S_0=\{(x_1,x_2)\in \mathcal S\mid\widetilde\pi(x_1)=\widetilde\pi(x_2)\}$ and we embed $L^{\infty}(\widetilde Y)\subset L^{\infty}(X)$ by $f\mapsto f\circ\widetilde\pi$.
\end{enumerate}

\end{proposition}

{\it Proof.} Since $\mathcal S$ has the Jones-Schmidt property, we can find a hyperfinite ergodic p.m.p. equivalence relation $\mathcal T$ on a probability space $(Y,\nu)$ and a factor map $\pi:(X,\mu)\rightarrow (Y,\nu)$ such that we have $\pi([x]_{\mathcal S})=[\pi(x)]_{\mathcal T}$, for almost every $x\in X$, and
$\mathcal S_0=\{(x_1,x_2)\in \mathcal S\mid\pi(x_1)=\pi(x_2)\}$ is strongly ergodic on almost every ergodic component of $\mathcal S_0$.

 We may assume that $\pi([x]_{\mathcal S})=[\pi(x)]_{\mathcal T}$, for all $x\in X$. We identify $L^{\infty}(X)^{[\mathcal S_0]}=L^{\infty}(\widetilde Y)$, for a probability space $(\widetilde Y,\tilde\nu)$. 
  Let $\widetilde\pi:(X,\mu)\rightarrow (\widetilde Y,\tilde\nu)$ and $\rho:(\widetilde Y,\tilde\nu)\rightarrow (Y,\nu)$ be the factor maps given by the embeddings $L^{\infty}(\widetilde Y)\subset L^{\infty}(X)$ and $L^{\infty}(Y)\subset L^{\infty}(\widetilde Y)$. Then  $\rho\circ\widetilde\pi=\pi$ and $\widetilde\pi(x_1)=\widetilde\pi(x_2)$, for every $(x_1,x_2)\in\mathcal S_0$. Therefore, $\mathcal S_0=\{(x_1,x_2)\in\mathcal S|\;\widetilde\pi(x_1)=\widetilde\pi(x_2)\}$. 

Consider the disintegration $\mu=\int_{\widetilde Y}\mu_y\;\text{d}\tilde\nu(y)$ of $\mu$ with respect to $\widetilde\pi$, where $\mu_y$ is a probability measure on $X$ supported on $\widetilde\pi^{-1}(\{y\})$, for every $y\in\widetilde Y$. 
  We claim that (b) holds, i.e., we have
  
  {\bf Claim 1.}  $(L^{\infty}(X)^{\omega})^{[\mathcal S_0]}=L^{\infty}(\widetilde Y)^{\omega}.$
   
{\it Proof of Claim 1.}
In order to prove the claim, it suffices to show that if  a sequence $f_n\in L^{\infty}(X)_1$ satisfies $\|\gamma f_n-f_n\|_2\rightarrow 0$, for every $\gamma\in [\mathcal S_0]$, then $\|f_n-E_{L^{\infty}(\widetilde Y)}(f_n)\|_2\rightarrow 0$. Let $\{\gamma_m\}_{m\geq 1}\subset [\mathcal S_0]$ be a sequence which is dense with respect to the metric $d(\gamma,\gamma')=\mu(\{x\in X\mid\gamma(x)\not=\gamma'(x)\})$. If the assertion concerning $(f_n)$ is false, then we can find a subsequence $(g_k)$ of $(f_n)$ such that $\|\gamma_mg_k-g_k\|_2\leq 1/2^k$, for all $k\geq m$, and $\inf_k\|g_k-E_{L^{\infty}(\widetilde Y)}(g_k)\|_2>0$. For $k\geq 1$ and $y\in\widetilde Y$, denote by $g_k^y\in L^{\infty}(\widetilde\pi^{-1}(\{y\}),\mu_y)_1$ the restriction of $g_k$ to $\pi^{-1}(\{y\})$. Then for every $m\geq 1$, we have that $$\int_{\widetilde Y}\sum_{k\geq1}\|\gamma_mg_k^y-g_k^y\|_{2,\mu_y}^2\;\text{d}\tilde\nu(y)=\sum_{k\geq 1}\|\gamma_mg_k-g_k\|_2^2<\infty.$$
Thus, for $\tilde\nu$-almost every $y\in\widetilde Y$, we have that $\|\gamma_mg_k^y-g_k^y\|_{2,\mu_y}\rightarrow 0$, for every $m\geq 1$.
Therefore, for $\tilde\nu$-almost every $y\in\widetilde Y$, we have that $\|\gamma g_k^y-g_k^y\|_{2,\mu_y}\rightarrow 0$, for every $\gamma\in [\mathcal S_0]$. Since $\mathcal S_0$ is $\mu_y$-strongly ergodic, we conclude that $\|g_k^y-\int g_k^y\;\text{d}\mu_y\|_{2,\mu_y}\rightarrow 0$, for $\tilde\nu$-almost every $y\in\widetilde Y$. From this we deduce that $\|g_k-E_{L^{\infty}(\tilde Y)}(g_k)\|_2^2=\int_{\widetilde Y}\|g_k^y-\int g_k^y\;\text{d}\mu_y\|_{2,\mu_y}\;\text{d}\nu(y)\rightarrow 0$, which gives a contradiction. \hfill$\square$

We continue with the following claim:

{\bf Claim 2.}
Let $\alpha:A\rightarrow B$ be an element of $[[\mathcal S]]$ such that there exists $\beta,\gamma\in[[\mathcal T]]$ satisfying $\pi\circ\alpha=\beta\circ\pi$ and $\pi\circ\alpha^{-1}=\gamma\circ\pi$.  Let $\widetilde A=\{y\in\widetilde Y|\mu_y(A)>0\}$ and $\widetilde B=\{y\in\widetilde Y|\mu_y(B)>0\}$. Then there exists a measure space isomorphism $\widetilde\alpha:(\widetilde A,{\widetilde\nu}_{|\widetilde A})\rightarrow (\widetilde B,{\widetilde\nu}_{|\widetilde B})$ such that $\widetilde\pi(\alpha(x))=\widetilde\alpha(\widetilde\pi(x))$, for almost every $x\in A$, and $\widetilde\pi(\alpha^{-1}(x))=\widetilde\alpha^{-1}(\widetilde\pi(x))$, for almost every $x\in B$.

{\it Proof of Claim 2.}
Note first that  \begin{equation}\label{disint}\int_{\widetilde Y}{\alpha_*}({\mu_y}_{|A})\;\text{d}\tilde\nu(y)={\alpha_*}(\mu_{|A})=\mu_{|B}=\int_{\widetilde Y}{\mu_y}_{|B}\;\text{d}\tilde\nu(y).\end{equation}
Let $y\in\widetilde A$ and $x_1,x_2\in\widetilde\pi^{-1}(\{y\})\cap A$.  Then $\widetilde\pi(x_1)=\widetilde\pi(x_2)$, hence $\pi(x_1)=\pi(x_2)$. Since $x_1,x_2\in A$, we get that $\pi(\alpha(x_1))=\pi(\alpha(x_2))$, thus $(\alpha(x_1),\alpha(x_2))\in\mathcal S_0$, and hence $\widetilde\pi(\alpha(x_1))=\widetilde\pi(\alpha(x_2))$. This proves that $\widetilde\pi$ is constant on $\alpha(\widetilde\pi^{-1}\{y\}\cap A)$. Thus, we can define a measurable map $\widetilde\alpha:\widetilde A\rightarrow \widetilde Y$ such that $\widetilde\pi(\alpha(x))=\widetilde\alpha(y)$, for all $y\in\widetilde A$ and $x\in \widetilde\pi^{-1}\{y\}\cap A$. This implies that $\widetilde\pi(\alpha(x))=\widetilde\alpha(\widetilde\pi(x))$, for almost every $x\in A$. Since $\alpha_*({\mu_y}_{|A})$ is supported on $\alpha(\widetilde\pi^{-1}(\{y\})\cap A)\subset\widetilde\pi^{-1}(\{\widetilde\alpha(y)\})$, by using \eqref{disint} and the uniqueness of the disintegration of $\mu_{|B}$ we conclude that ${\mu_{\widetilde\alpha(y)}}_{|B}=\alpha_*({\mu_y}_{|A})$ and $\widetilde\alpha(y)\in\widetilde B$, for almost every $y\in\widetilde A$, and that $\widetilde\alpha:\widetilde A\rightarrow\widetilde B$ is onto.
Similarly, we can find an onto measurable map $\overline{\alpha}:\widetilde B\rightarrow\widetilde A$ such that $\widetilde\pi(\alpha^{-1}(x))=\overline{\alpha}(\widetilde\pi(x))$, for almost every $x\in B$. It is easy to see that $\overline{\alpha}=\widetilde\alpha^{-1}$.

If $C\subset\widetilde B$ is a measurable set, then $\widetilde\pi^{-1}(\widetilde{\alpha}^{-1}(C))=\{x\in A|\widetilde\alpha(\widetilde\pi(x))\in C\}=\{x\in A|\widetilde\pi(\alpha(x))\in C\}$. Since $\alpha$ is $\mu$-preserving, we get that $$\widetilde\nu(\widetilde\alpha^{-1}(C))=\mu(\widetilde\pi^{-1}(\widetilde\alpha^{-1}(C)))=\mu(\{x\in A|\widetilde\pi(\alpha(x))\in C\})=\mu(\{x\in B|\widetilde\pi(x)\in C\})=\widetilde\nu(C).$$ This shows that $\widetilde\alpha$ is  $\widetilde\nu$-preserving, and finishes the proof of the claim. \hfill$\square$

Now, the first assumption on $\pi$ implies we can find a sequence of elements $\alpha_n:A_n\rightarrow B_n, n\geq 1$, of $[[\mathcal S]]$ which satisfy the hypothesis of Claim 2 and such that $[x]_{\mathcal S}=\{\alpha_n(x)\}_{n\geq 1}$, for almost every $x\in X$. For $n\geq 1$, let $\widetilde\alpha_n:\widetilde A_n\rightarrow\widetilde B_n$ be the $\widetilde\nu$-preserving isomorphism obtained by applying Claim 2 to $\alpha_n$. 

Let $\widetilde{\mathcal T}$ be the countable p.m.p. equivalence relation of 
$(\widetilde Y,\widetilde\nu)$ generated by $\{\widetilde\alpha_n\}_{n\geq 1}$.
It is then immediate that condition (a) is satisfied. 
If $g\in L^{\infty}(\widetilde Y)$ is $\widetilde{\mathcal T}$-invariant, then $g\circ\widetilde\pi$ is $\mathcal S$-invariant. Since $\mathcal S$ is ergodic, we get that $g\circ\widetilde\pi\in\mathbb C1$, hence $g\in\mathbb C1$, showing that $\widetilde{\mathcal T}$ is ergodic.

Finally, we show that $\widetilde{\mathcal T}$ is hyperfinite. 
Since $\mathcal T$ is hyperfinite, it suffices to show that  the restriction of $\rho$ to $[y]_{\widetilde{\mathcal T}}$ is injective and $\rho([y]_{\widetilde{\mathcal T}})\subset [\rho(y)]_{\mathcal T}$, for all $y\in\widetilde Y$. To this end, let $y\in\widetilde Y$ and $y'\in [y]_{\widetilde{\mathcal T}}$. Then we can find $(x,x')\in\mathcal S$ such that $y=\widetilde\pi(x)$ and $y'=\widetilde\pi(y)$. Since $\rho(y)=\pi(x)$, $\rho(y')=\pi(x')$, and  $\pi(x')\in [\pi(x)]_{\mathcal T}$, we get that $\rho(y')\in[\rho(y)]_{\mathcal T}$.  If $\rho(y')=\rho(y)$, then $\pi(x')=\pi(x)$, hence $(x,x')\in\mathcal S_0$ thus $y=\widetilde\pi(x)=\widetilde\pi(x')=y'$, showing that $\rho_{|[y]_{\widetilde{\mathcal T}}}$ is injective. This finishes the proof. \hfill$\blacksquare$


\subsection*{Proof of Theorem \ref{D}} Since $N'\cap N^{\omega}\subset A^{\omega}$, Theorem \ref{A} implies the existence of $u\in\mathcal U(P\bar{\otimes}R)$,  $t>0$, and a Cartan subalgebra $B\subset P^t$ such that $\theta(A)=uBu^*$, where we identify $P\bar{\otimes}R=P^t\bar{\otimes}R^{1/t}$. After replacing $\theta$ with  $\text{Ad}(u^*)\circ\theta$, we may assume that $\theta(A)=B$.
We also identify $P^t\bar{\otimes}R^{1/t}=P^t\bar{\otimes}R$ using an isomorphism $R^{1/t}\mapsto R$. Thus, we see $\theta$ as an isomorphism $\theta:N\bar{\otimes}R\rightarrow P^t\bar{\otimes}R$ satisfying $\theta(A)=B$. 
Finally, we identify $A=L^{\infty}(X), B=L^{\infty}(Z)$, for some probability spaces $(X,\mu), (Z,\eta)$. Let $\alpha:(X,\mu)\rightarrow (Z,\eta)$ be the measure space isomorphism given by $\theta(a)=a\circ\alpha^{-1}$, for every $a\in A$.

Denote $\mathcal S=\mathcal R(A\subset N)$. For $g\in[\mathcal S]$,  let $u_g\in\mathcal N_{N}(A)$ such that $u_gau_g^*=\sigma_g(a):=a\circ g^{-1}$, for every $a\in A$.
By Proposition \ref{cocycle}, $\theta(A\bar{\otimes}R)=B\bar{\otimes}R$ and we can find
 $w_g\in \mathcal U(A\bar{\otimes}R)$, for every $g\in [\mathcal S]$, and a measurable map $\varphi:X\rightarrow\text{Aut}(R)$ such that the following conditions hold:
\begin{enumerate}
\item $\theta(w_gu_g)\in P^t$, for all $g\in [\mathcal S]$.
\item $\theta^{-1}(T)(x)=\varphi_{x}(T(\alpha(x))$, for every $T\in B\bar{\otimes}R$ and almost every $x\in X$.
\item $\text{Ad}(w_g(x))=\varphi_x\varphi_{g^{-1}x}^{-1}$, for all $g\in [\mathcal S]$ and almost every $x\in X$.
\end{enumerate}
Since $\mathcal S$ has the Jones-Schmidt property, Proposition \ref{eqJS} provides a hyperfinite ergodic p.m.p. equivalence relation $\mathcal T$ on a probability space $(Y,\nu)$ and a factor map $\pi:(X,\mu)\rightarrow (Y,\nu)$ such that 
 $\pi([x]_{\mathcal S})=[\pi(x)]_{\mathcal T}$, for almost every $x\in X$, and if we let $\mathcal S_0=\{(x_1,x_2)\in \mathcal S\mid\pi(x_1)=\pi(x_2)\}$, then $(L^{\infty}(X)^{\omega})^{[\mathcal S_0]}=L^{\infty}(Y)^{\omega}$. In particular, $L^{\infty}(X)^{[\mathcal S_0]}\subset L^{\infty}(Y)$ and 
 $(L^{\infty}(X)^{\omega})^{[\mathcal S]}\subset L^{\infty}(Y)^{\omega}$.

Note that (3) implies that $\varphi_x\varphi_y^{-1}\in\text{Inn}(R)$, for almost every $(x,y)\in\mathcal S$. By applying Theorem \ref{rigid}, we get measurable maps $\psi:Y\rightarrow\text{Aut}(R)$ and $U:X\rightarrow\mathcal U(R)$ such that \begin{equation}\label{psi}\varphi_x=\text{Ad}(U(x))\psi_{\pi(x),}\;\;\;\;\text{for almost every $x\in X$}.\end{equation} 

View the function $U:X\rightarrow \mathcal U(R)$ as an element of $\mathcal U(A\bar{\otimes}R)$. 
Let $\Theta:N\bar{\otimes}R\rightarrow P^t\bar{\otimes}R$ be the isomorphism given by $\Theta=\theta\circ\text{Ad}(U)$. For $g\in [\mathcal S]$, let $W_g=U^*w_g(\sigma_g\otimes\text{id}_R)(U)\in\mathcal U(A\bar{\otimes}R)$. By using \eqref{psi}, the above conditions (1)-(3) rewrite as:
\begin{enumerate}[label=(\alph*)]
\item\label{a} $\Theta(W_gu_g)\in P^t$, for all $g\in [\mathcal S]$.
\item $\Theta^{-1}(T)(x)=\psi_{\pi(x)}(T(\alpha(x))$, for every $T\in B\bar{\otimes}R$ and almost every $x\in X$.
\item $\text{Ad}(W_g(x))=\psi_{\pi(x)}\psi_{\pi(g^{-1}x)}^{-1}$, for all $g\in [\mathcal S]$ and almost every $x\in X$.
\end{enumerate}

By using condition (c) we get that $\psi_y\psi_{z}^{-1}\in\text{Inn}(R)$, for  almost every $(y,z)\in \mathcal T$.
We claim that there is a cocycle $c:\mathcal T\rightarrow\mathcal U(R)$ such that  \begin{equation}\label{psi3}\psi_y\psi_z^{-1}=\text{Ad}(c(y,z)),\;\;\;\text{for almost every $(y,z)\in\mathcal T$.}\end{equation}
Since $\mathcal T$ is ergodic hyperfinite, we can find a free ergodic p.m.p. action $\mathbb Z\curvearrowright (Y,\nu)$ such that $\mathcal T=\mathcal R(\mathbb Z\curvearrowright Y)$. 
Let $T_0:Y\rightarrow Y$ be the generator of the action $\mathbb Z\curvearrowright Y$.
Since $T_0\in [\mathcal T]$, by (c) we can find a measurable map $d:Y\rightarrow\mathcal U(R)$ such that $\psi_{T_0(y)}\psi_y^{-1}=\text{Ad}(d(y))$, for almost every $y\in Y$.  
Then the claim holds for the unique cocycle $c:\mathcal T\rightarrow\mathcal U(R)$ satisfying $c(T_0(y),y)=d(y)$, for all $y\in Y
$.

 For $g\in [\mathcal S]$, we define $\widetilde W_g\in A\bar{\otimes}R$ by letting $\widetilde W_g(x)=c(\pi(x),\pi(g^{-1}x))$. By combining (c) and \eqref{psi3} we get that $\text{Ad}(W_g(x))=\text{Ad}(\widetilde W_g(x))$, for almost every $x\in X$, and thus $\widetilde W_gW_g^*\in A$.

The rest of proof relies on the following claim:

{\bf Claim.} There exists $\Psi\in\overline{\text{Inn}(N\bar{\otimes}R)}$ such that  $\Psi(A\bar{\otimes}R)=A\bar{\otimes}R$, $\Psi(T)(x)=\psi_{\pi(x)}(T(x))$, for all $T\in A\bar{\otimes}R$ and almost every $x\in X$, and $\Psi(u_g)=\widetilde W_gu_g$, for all $g\in [\mathcal S]$.

{\it Proof of the claim.}
Write $R=\bar{\otimes}_{k\in\mathbb N}\mathbb M_2(\mathbb C)$ and define $R^n=\bar{\otimes}_{k<n}\mathbb M_2(\mathbb C)$, for $n\in\N$. Fix $n\in\N$ and consider the set $B_n:=\{(y,u)\in Y\times\cU(R)\mid {\psi_y}_{|R^n} = \text{Ad}(u)_{|R^n}\}$. As $\psi:Y\rightarrow\Aut(R)$ is measurable, we see that $B_n\subset Y\times\cU(R)$ is a Borel subset. Moreover we have that $\pi_1(B_n)=Y$, where $\pi_1$ is the projection onto the first component. Applying \cite[Theorem~A.16]{Ta01} we get a measurable map $u_n:Y\rightarrow\mathcal U(R)$ such that ${\psi_y}_{|R^n}=\text{Ad}(u_n(y))_{|R^n}$, for all $y\in Y$. Thus, $\lim\limits_{n\rightarrow\infty}\text{Ad}(u_n(y))(v)=\psi_y(v)$, for all $y\in Y$ and $v\in R$.

Define a cocycle $c_n:\mathcal T\rightarrow\mathcal U(R)$ by $c_n(y,z)=u_n(y)^*c(y,z)u_n(z)$. Then $\text{Ad}(c_n(y,z))_{|R^n}=\text{id}_{R^n}$, for all $(y,z)\in\mathcal T$. Since $\mathcal T$ is hyperfinite, by applying Lemma \ref{trivial}, we can find a subsequence $\{n_k\}_{k\geq 1}$ of $\mathbb N$ and measurable maps $d_k:Y\rightarrow\mathcal U(R)$, for $k\geq 1$, such that $\|c_{n_k}(y,z)-d_k(y)d_k(z)^{-1}\|_2\rightarrow 0$ and $\|\text{Ad}(d_k(y))(v)-v\|_2\rightarrow 0$, as $k\rightarrow\infty$, for almost every $(y,z)\in\mathcal T$ and every $v\in R$.

For $k\geq 1$, define $U_k\in\mathcal U(A\bar{\otimes}R)$ by letting $U_k(x)=u_{n_k}(\pi(x))d_k(\pi(x))$.
Then for every $g\in [\mathcal S]$ and almost every $x\in X$ we have that 
\begin{equation} (U_k(\sigma_g\otimes\text{id}_R)(U_k^*))(x)=u_{n_k}(\pi(x))d_k(\pi(x))d_k(\pi(g^{-1}x))^*u_{n_k}(\pi(g^{-1}x))^*\rightarrow \widetilde W_g(x). \end{equation}

Thus, $\|U_ku_gU_k^*-\widetilde W_gu_g\|_2\rightarrow 0$, for all $g\in [\mathcal S]$. 
If $T\in R$, then as $\|\text{Ad}(d_k(\pi(x))(T)-T\|_2\rightarrow 0$, for almost every $x\in X$, we get that $\lim\limits_{k\rightarrow\infty}\text{Ad}(U_k(x))(T)=\lim\limits_{k\rightarrow\infty}\text{Ad}(u_{n_k}(\pi(x)))(T)=\psi_{\pi(x)}(T).$  Also, $U_kTU_k^*=T$, for all $T\in A$. Since $\{u_g\}_{g\in [\mathcal S]}\cup R\cup A$ generates $N\bar{\otimes}R$ as a von Neumann algebra we conclude that the limit $\Psi=\lim\limits_{k\rightarrow\infty}\text{Ad}(U_k)\in\overline{\text{Inn}(N\bar{\otimes}R)}$ exists and satisfies the claim. 
\hfill$\square$

Finally, let $g\in [\mathcal S]$. Since $\widetilde W_gW_g^*\in A$, we have $\Theta(\widetilde W_gW_g^*)=\theta(\widetilde W_gW_g^*)\in B$. Since $\Theta(W_gu_g)\in P^t$ by (a), we get that  $(\Theta\circ\Psi)(u_g)=\Theta(\widetilde W_gu_g)=\Theta(\widetilde W_gW_g^*)\Theta(W_gu_g)\in P^t.$ 
Since $(\Theta\circ\Psi)(A)=B$ and $\{u_g\}_{g\in[\mathcal S]}\cup A$ generates $N$ as a von Neumann algebra, we deduce that $(\Theta\circ\Psi)(N)\subset P^t$.
Let $T\in A\bar{\otimes}R$ and define $\widetilde T\in B\bar{\otimes}R$ by letting $\widetilde T(z)=T(\alpha^{-1}(z))$, for all $z\in Z$. Then by combining (b) and the claim we get that $\Theta^{-1}(\widetilde T)(x)=\psi_{\pi(x)}(\widetilde T(\alpha(x))=\psi_{\pi(x)}(T(x))=\Psi(T)(x)$, for almost every $x\in X$. This shows that $(\Theta\circ\Psi)(T)=T\circ\alpha^{-1},$ for every $T\in A\bar{\otimes}R$, and in particular implies that $(\Theta\circ\Psi)(R)=R$. 
Since $\Theta\circ\Psi:N\bar{\otimes}R\rightarrow P^t\bar{\otimes}R$ is onto and $(\Theta\circ\Psi)(N)\subset P^t$, we derive that $(\Theta\circ\Psi)(N)=P^t$.
Thus, we have that $\Theta\circ\Psi=\theta_1\otimes\theta_2$, where $\theta_1:N\rightarrow P^t$ and $\theta_2:R\rightarrow R$ are isomorphisms.
This implies the conclusion of Theorem \ref{D} for $s=1/t$.
\hfill$\blacksquare$

\section{Equivalence relations without the Jones-Schmidt property}

We begin with the following lemma, which will be key in the proofs of Theorems \ref{E} and \ref{E'} below.
\begin{lemma}\label{aux}
Let $\mathcal S$ be a countable ergodic p.m.p. equivalence relation on a probability space $(X,\mu)$ with the Jones-Schmidt property. Denote $M=L(\mathcal S)$ and $A=L^{\infty}(X)$. Assume that $A_n\subset A$ are von Neumann subalgebras such that $A=\bar{\otimes}_{n\in\mathbb N}A_n$ and $M'\cap A^{\omega}\subset (\bar{\otimes}_{k\geq n}A_k)^{\omega}$, for every $n\in\mathbb N$.

Then for every $\varepsilon>0$, we can find finite dimensional subalgebras $\tilde A_n\subset A_n$, for $n\in\mathbb N$, such that $$M'\cap A^{\omega}\subset_{\varepsilon}(\bar{\otimes}_{n\in\mathbb N}\tilde A_n)^\omega.$$
\end{lemma}

{\it Proof.} Let $(Y,\nu)=(Y_0,\nu_0)^{\mathbb N}$, where $Y_0=\{0,1\}$ and $\nu_0=\frac{1}{2}(\delta_0+\delta_1)$, and consider the equivalence relation $\mathcal T$ on $(Y,\nu)$ given by $(y_k)\mathcal T(z_k)$ if and only if $y_k=z_k$, for all but finitely many $k\in\mathbb N$.
Since $\mathcal S$ has the Jones-Schmidt property, by applying Proposition \ref{eqJS} we can find a factor map $\pi:(X,\mu)\rightarrow (Y,\nu)$ such that  
(a) $\pi([x]_{\mathcal S})=[\pi(x)]_{{\mathcal T}}$, for almost every $x\in X$, and
(b) $(A^{\omega})^{[\mathcal S_0]}=B^{\omega}$, where we embed $B=L^{\infty}(Y)$ into $A$ via the map $f\mapsto f\circ\pi$ and denote  $\mathcal S_0=\{(x_1,x_2)\in \mathcal S\mid\pi(x_1)=\pi(x_2)\}$.

Note first that since $M'\cap A^{\omega}\subset (A^{\omega})^{[\mathcal S]}\subset (A^{\omega})^{[\mathcal S_0]}$,  condition (b) implies that
\begin{equation}\label{centr}
M'\cap A^{\omega}\subset B^{\omega}.
\end{equation}
 Next, we identify $B=\bar{\otimes}_{k\in\mathbb N}L^{\infty}(Y_0,\nu_0)_k$ and put $B_n=\bar{\otimes}_{k\geq n}L^{\infty}(Y_0,\nu_0)_k$, for $n\in\mathbb N$. 
 We claim that \begin{equation}\label{centr1} \prod_{\omega}B_n\subset M'\cap A^{\omega}.\end{equation}

To see this,
let $f\in \prod_{\omega}B_n$ and $\alpha\in [\mathcal S]$.
If we view $f$ as an element of $A^{\omega}$, then we can represent it as $f=(f_n\circ\pi)$, where $f_n\in B_n$, for every $n\in\mathbb N$, and $\sup\|f_n\|<\infty$.  By  condition (a) above, we can find a measurable map $\beta:Y\rightarrow Y$ such that  $\pi(\alpha^{-1}(x))=\beta(\pi(x))$ and $\beta(y)\in [y]_{\mathcal T}$, for almost every $x\in X$ and $y\in Y$. 
Fix $n\in\mathbb N$, and denote $Y_n=\{\text{$(y_k)\in Y|\beta(y)_k=y_k$, for all $k\geq n$}\}$. Then we have that $f_n(\beta(y))=f_n(y)$, for all  $y\in Y_n$. 
Since $u_{\alpha}(f_n\circ\pi)u_{\alpha}^*=f_n\circ\pi\circ{\alpha}^{-1}=f_n\circ\beta\circ\pi$, we get that
$$\|u_{\alpha}(f_n\circ\pi)u_{\alpha}^*-f_n\circ\pi\|_2=\|f_n\circ\beta\circ\pi-f_n\circ\pi\|_2=\|f_n\circ\beta-f_n\|_2\leq 2\sqrt{\nu(Y\setminus Y_n)}\|f_n\|.$$

Since $\lim\limits_{n\rightarrow\infty}\nu(Y_n)=1$, we conclude that $\lim\limits_{n\rightarrow\infty}\|u_{\alpha}(f_n\circ\pi)u_{\alpha}^*-f_n\circ\pi\|_2=0$, and hence $u_{\alpha}fu_{\alpha}^*=f$. Since this holds for every $f\in\prod_{\omega}B_n$ and $\alpha\in [\mathcal S]$, claim \eqref{centr1} follows.

For $n\in\mathbb N$, denote $C_n=\bar{\otimes}_{k\geq n}A_k$. We claim that \begin{equation}\label{strongly} \text{$B\prec_{A}^sC_n$, for every $n\in\mathbb N$.} \end{equation}
To see this, let $n\in\mathbb N$ and $p\in A$ be a non-zero projection. Put $\delta=\uptau(p)/2>0$.
Since $M'\cap A^{\omega}\subset C_n^{\omega}$, by combining \eqref{centr1} with Lemma \ref{easy}, we deduce that $B_m\subset_{\delta}C_n$, for some $m\in\mathbb N$. Lemma \ref{epsilon} then implies the existence of a projection $q\in A$ such that $B_mq\prec_{A}^sC_n$ and $\uptau(q)\geq 1-\delta$. 
Since $B=(\bar{\otimes}_{k<m}L^{\infty}(Y_0,\nu_0)_k)\bar{\otimes}B_m$, it follows that $Bq\prec_{A}^sC_{n}$. Since $\uptau(p)+\uptau(q)>1$, we have that $r=pq\in Ap$ is a non-zero projection such that $Br\prec_{A}C_n$. This proves the claim made in \eqref{strongly}.

Now, for $n\in\mathbb N$, let $\{A_{n,k}\}_{k\in\mathbb N}$ be an increasing sequence of finite dimensional subalgebras of $A_n$ such that $\cup_{k\in\mathbb N}A_{n,k}$ is weakly dense in $A_n$. Fix $\varepsilon>0$ and $n\in\mathbb N$. Since $A$ is abelian, \eqref{strongly} implies the existence 
of projections $\{p_{n,l}\}_{l\in\mathbb N}$ in $A$ such that $\sum_{l\in\mathbb N}p_{n,l}=1$ and $(B)_1p_{n,l}\subset (C_n)_1p_{n,l}$, for all $l\in\mathbb N$.
We claim that we can find $k_{n,1},k_{n,2},...,k_{n,n-1}\in \mathbb N$ such that \begin{equation}\label{containment} \text{$B\subset_{\varepsilon/2^n}\mathcal A_n:=(\bar\otimes_{1\leq m\leq n-1}A_{m,k_{n,m}})\bar{\otimes}C_n$, for every $n\in\mathbb N$.} \end{equation}
To justify this claim, let $N\in\mathbb N$ such that $\|\sum_{l>N}p_{n,l}\|_2<\varepsilon/2^{n+1}$. Since $A=(\bar{\otimes}_{1\leq m\leq n-1}A_m)\bar{\otimes}C_n$, we have that $\cup_{k_1,...,k_{n-1}\in\mathbb N}\big((\bar{\otimes}_{1\leq m\leq n-1}A_{m,k_m})\bar{\otimes}C_n\big)$ is weakly dense in $A$. Thus, we can find $k_{n,1},k_{n,2},...,k_{n,n-1}\in \mathbb N$ such that  $\|p_{n,l}-E_{\mathcal A_n}(p_{n,l})\|_2\leq \varepsilon/(2^{n+1}N)$, for all $1\leq l\leq N$, where $\mathcal A_n$ is defined as in \eqref{containment}. Let $x\in (B)_1$. Then for every $1\leq l\leq N$, we can find $y_l\in (C_n)_1$ such that $xp_{n,l}=y_lp_{n,l}$. Since $\sum_{l=1}^Ny_lE_{\mathcal A_n}(p_{n,l})\in\mathcal A_n$, \eqref{containment} is a consequence of the following estimate: \begin{align*}
\|x-\sum_{l=1}^Ny_lE_{\mathcal A_n}(p_{n,l})\|_2\leq \|x-\sum_{l=1}^Ny_lp_{n,l}\|_2+\varepsilon/2^{n+1}=\|x-\sum_{l=1}^Nxp_{n,l}\|_2+\varepsilon/2^{n+1}\leq\varepsilon/2^n.
\end{align*}

Put $\mathcal A=\cap_{n\in\mathbb N}\mathcal A_n$ and let $x\in (B)_1$.  Then \eqref{containment} gives that $\|x-E_{\mathcal A_k}(x)\|_2\leq\varepsilon/2^k$, for all $k\in\mathbb N$. Since $E_{A_1}\circ E_{A_2}\circ...\circ E_{A_n}=E_{\cap_{k=1}^nA_k}$, the last inequality implies that $\|x-E_{\cap_{k=1}^n\mathcal A_k}(x)\|_2\leq\varepsilon$, for all $n\in\mathbb N$.
Since $\lim\limits_{n\rightarrow\infty}\|E_{\cap_{k=1}^n\mathcal A_n}(x)-E_{\mathcal A}(x)\|_2=0$, we conclude that $\|x-E_{\mathcal A}(x)\|_2\leq\varepsilon$. 
As this holds for all $x\in (B)_1$, we deduce that $B\subset_{\varepsilon}\mathcal A$. In combination with \eqref{centr}, it follows that $M'\cap A^{\omega}\subset_{\varepsilon}\mathcal A^{\omega}$.
Finally, note that $\mathcal A=\bar{\otimes}_{m\in\mathbb N}A_{m,k_m}$, where $k_m=\min\{k_{n,m}|n\geq m+1\}$, for every $m\in\mathbb N$.
This implies the conclusion of the lemma.
\hfill$\blacksquare$

\subsection*{Proof of Theorem \ref{E}} Assume by contradiction that $\mathcal S:=\mathcal R(\Gamma\curvearrowright X)$ has the Jones-Schmidt property. Denote $A=L^{\infty}(X)$ and $M=A\rtimes\Gamma$.
For $n\in\mathbb N$, let $A_n=L^{\infty}(X_n)$ and $M_n=A_n\rtimes\Gamma_n$. Then $A=\bar{\otimes}_{n\in\mathbb N}A_n$. For $n\in\mathbb N$, let $C_n=\bar{\otimes}_{k\geq n}A_k$ and $D_n=\bar{\otimes}_{k\not=n}A_k$. We claim that \begin{equation}\label{centr2} \text{$M'\cap A^{\omega}\subset C_n^{\omega}$.}
 \end{equation}
 Indeed, if $k\in\mathbb N$, then since the action $\Gamma_k\curvearrowright (X_k,\mu_k)$ is strongly ergodic we have that $M_k'\cap A_k^{\omega}=\mathbb C1$. This easily implies that $M'\cap A^{\omega}\subset M_k'\cap A^{\omega}\subset D_k^{\omega}$ (by using, e.g., \cite[Proposition 5.2]{Ma17}). Thus, for every $n\in\mathbb N$, we have that $M'\cap A^{\omega}\subset\cap_{k=1}^{n-1}D_k^{\omega}=C_n^{\omega}$ and claim \eqref{centr2} follows.
 
We may thus apply Lemma \ref{aux} to deduce the existence of finite dimensional subalgebras $\tilde A_n\subset A_n$, for $n\in\mathbb N$, such that $M'\cap A^{\omega}\subset_{1/2}\mathcal A^{\omega}$, where $\mathcal A=\bar{\otimes}_{n\in\mathbb N}\tilde A_n$. On the other hand, since $\tilde A_n$ is finite dimensional and $A_n$ is diffuse, we can find $a_n\in\mathcal U(A_n)$ such that $E_{\tilde A_n}(a_n)=0$, for all $n\in\mathbb N$. 
Let $u_n\in\mathcal U(A)$ be given by $u_n=\otimes_{k\in\mathbb N}u_{n,k}$, where $u_{n,n}=a_n$ and $u_{n,k}=1$, if $k\not=n$.
 Since $E_{\tilde A_n}(a_n)=0$, we get that $E_{\mathcal A}(u_n)=0$. Hence, if we put $u=(u_n)_n\in \mathcal U(A^{\omega})$, then $E_{\mathcal A^{\omega}}(u)=0$. Since we clearly have that $u\in M'\cap A^{\omega}$ this contradicts the fact that $M'\cap A^{\omega}\subset_{1/2}\mathcal A^{\omega}$.
 \hfill$\blacksquare$
 
 \subsection*{Proof of Theorem \ref{E'}} Assume by contradiction that $\mathcal S:=\mathcal R(\Gamma\curvearrowright X)$ has the Jones-Schmidt property. 
Denote $A=L^{\infty}(X)$ and $M=A\rtimes\Gamma$.
For $n\in\mathbb N$, let $A_n=L^{\infty}(X_n)$. Then $A=\bar{\otimes}_{n\in\mathbb N}A_n$. 
Let $F_{n,k}\in (A_n)_1$,  $n,k\in\mathbb N$, be as in the hypothesis of Theorem \ref{E'} and put $c:=\inf_{n,k}\|F_{n,k}\|_2>0$.

For $n\in\mathbb N$, let $C_n=\bar{\otimes}_{k\geq n}A_k$ and $D_n=\bar{\otimes}_{k\not=n}A_k$. We claim that \begin{equation}\label{centr3} \text{$M'\cap A^{\omega}\subset C_n^{\omega}$.}
 \end{equation}
 To see this, let $n\in\mathbb N$. Denote by $\pi_n$ the Koopman representation of $\Gamma$ on $L^2(X_n)\ominus\mathbb C1$.  Then $\pi_n\otimes\bar{\pi}_n$ is a subrepresentation of the Koopman representation of $\Gamma$ on $L^2(X_n\times X_n)\ominus\mathbb C1$. Since the latter representation is assumed to have spectral gap, we deduce that $\pi_n\otimes\bar{\pi}_n$ has spectral gap. By applying \cite[Lemmas 3.2 and 3.5]{Po06a}, we get that $M'\cap A^{\omega}\subset D_n^{\omega}$. Thus, for every $n\in\mathbb N$, we have that $M'\cap A^{\omega}\subset\cap_{k=1}^{n-1}D_k^{\omega}=C_n^{\omega}$ and claim \eqref{centr3} follows.
 
We may thus apply Lemma \ref{aux} to deduce the existence of finite dimensional subalgebras $\tilde A_n\subset A_n$, for $n\in\mathbb N$, such that $M'\cap A^{\omega}\subset_{c/3}\mathcal A^{\omega}$, where $\mathcal A=\bar{\otimes}_{n\in\mathbb N}\tilde A_n$.  We claim that there is $n\in\mathbb N$ such that \begin{equation}\label{centr4}\text{$\|F_{n,k}-E_{\tilde A_n}(F_{n,k})\|_2\leq c/2$, for every $k\in\mathbb N$.}\end{equation} If this were false, then for every $n\in\mathbb N$ we can find $k_n\in\mathbb N$ such that $\|F_{n,k_n}-E_{\tilde A_n}(F_{n,k_n})\|_2>c/2$. But then $F=(F_{n,k_n})_n\in (A^{\omega})_1$ satisfies $\|F-E_{\cA^\omega}(F)\|_2\geq c/2$. On the other hand, the hypothesis implies that $F\in M'\cap A^{\omega}$, which contradicts the fact that $M'\cap A^{\omega}\subset_{c/3}\mathcal A^{\omega}$ and thus proves \eqref{centr4}.

Finally, since $F_{n,k}\rightarrow 0$ weakly in $L^2(X_n)$ and $\tilde A_n$ is finite dimensional, we get that $\|E_{\tilde A_n}(F_{n,k})\|_2\rightarrow 0$, as $k\rightarrow\infty$. Since $c>0$, this contradicts \eqref{centr4} and finishes the proof. \hfill$\blacksquare$

We end this section by justifying the claim made in Example \ref{free}.

\begin{remark}\label{freeproof} 
Recall that $\widetilde{\pi}_t=\oplus_{i\in\mathbb N}\pi_t:\Gamma\rightarrow\mathcal O(\widetilde{\mathcal H}_t)$, where $\pi_t$ is the orthogonal GNS representation of $\Gamma=\mathbb F_m$ associated to the positive definite function $\varphi_t(g)=e^{-t|g|}$, for $t>0$.  Let $\Gamma\curvearrowright (X_t,\mu_t)$ be the Gaussian representation associated to $\widetilde{\pi}_{t}$. Our goal is to show that the diagonal action $\Gamma\curvearrowright (X,\mu):=\prod_{n\in\mathbb N}(X_{t_n},\mu_{t_n})$ satisfies the hypothesis of Theorem \ref{E'}, for any sequence $t_n\rightarrow 0$. 

First, we claim that the diagonal action $\Gamma\curvearrowright X_t\times X_t$ has spectral gap, for every $t>0$. To this end, assume that $m<\infty$. Since the Koopman representation of $\Gamma$ on $L^2(X_t)\ominus\mathbb C1$  is isomorphic to a subrepresentation of a multiple of $\bigoplus_{k\in\mathbb N}\pi_t^{\otimes k}$, the same is true for the Koopman representation of $\Gamma$ on $L^2(X_t\times X_t)\ominus\mathbb C1$. Since $m<\infty$, we have that  $\varphi_t^{N}\in\ell^1(\Gamma)$, for some $N\geq 1$, hence, $\pi_t^{\otimes N}$ is contained in the left regular representation of $\Gamma$. Since $\Gamma$ is non-amenable, this implies that $\bigoplus_{k\in\mathbb N}\pi_t^{\otimes k}$ has spectral gap. In combination with the above it follows that the diagonal action $\Gamma\curvearrowright X_t\times X_t$ has spectral gap. If $m=+\infty$, then the same conclusion can be reached by replacing $\Gamma
$ with a finitely generated non-abelian subgroup in the above argument.

Secondly, let $t\in\mathbb N$. Note that we can find pairwise orthogonal vectors $\{\xi_{t,k}\}_{k\in\mathbb N}$ in $\tilde{\mathcal H}_{t}$ such that  $\langle\widetilde{\pi}_{t}(g)\xi_{t,k},\xi_{t,k}\rangle=\varphi_{t}(g)$, for all $g\in\Gamma$ and $k\in\mathbb N$. Recall that there is a map $u:\tilde{\mathcal H}_{t}\rightarrow\mathcal U(L^{\infty}(X_{t}))$ such that $u(\widetilde{\pi}_t(g)\xi)=u(\xi)\circ g^{-1}$, $u(\xi+\eta)=u(\xi)u(\eta)$, $\overline{u(\xi)}=u(-\xi)$, and $\int_{X_t}u(\xi)=\exp{(-\|\xi\|_2^2/2)}$, for all $g\in\Gamma$ and $\xi,\eta\in\widetilde{\mathcal H}_{t}$. For $k\in\mathbb N$, define $$F_{t,k}:=u(\xi_{t,k})-\int_{X_t}u(\xi_{t,k})\in L^{\infty}(X_t).$$
Then we have $\|F_{t,k}\|_{\infty}\leq 2$, $\|F_{t,k}\|_2=\sqrt{1-\exp{(-1)}}$, and $\|F_{t,k}\circ g-F_{t,k}\|_2=\sqrt{2-2\exp{(-1+\varphi_t(g))}}$, for all $k\in\mathbb N$ and $g\in\Gamma$. Moreover, $\int_{X_t}F_{t,k}\overline{F_{t,k'}}=0$, for all $k\not=k'$, hence $F_{t,k}\rightarrow 0$ weakly in $L^2(X_t)$. It is now clear that $F_{n,k}:=F_{t_n,k}/2\in L^{\infty}(X_{t_n})$ satisfy the hypothesis of Theorem \ref{E'}.
\end{remark}

\section{Equivalence relations with a unique stable decomposition}

We start off with the following easy lemma, which we will use in the proof of Theorem~\ref{F} below.

\begin{lemma}\label{lem:split}
	Suppose $M$, $N$ are II$_1$ factors, and $A\subset M$, $B\subset N$ are Cartan subalgebras. Then
	\[
	\cR(A\subset M)\times \cR(B\subset N) \cong \cR(A\otb B\subset M\otb N).
	\]
\end{lemma}

{\it Proof.} We identify $A$ with $L^\infty(X)$ and $B$ with $L^\infty(Y)$ for probability spaces $(X,\mu)$ and $(Y,\nu)$, and we write $\cR:=\cR(A\subset M)$ and $\cS:=\cR(B\subset N)$. Following \cite{FM77} we can find 2-cocycles $v\in$ H$^2(\cR,\mathbb T)$ and $w\in$ H$^2(\cS,\mathbb T)$ such that
\[
(A\subset M) \cong (L^\infty(X)\subset L_v(\cR)), \qquad\text{and}\qquad (B\subset N) \cong (L^\infty(Y)\subset L_w(\cS)).
\]
It follows that
\[
(A\otb B\subset M\otb N) \cong (L^\infty(X)\otb L^\infty(Y) \subset L_v(\cR)\otb L_w(\cS)) \cong (L^\infty(X\times Y)\subset L_{v\times w}(\cR\times \cS)),
\]
where $v\times w\in$ H$^2(\cR\times\cS,\mathbb T)$ is the cocycle defined by
\[
(v\times w)((x_1,y_1),(x_2,y_2),(x_3,y_3)) = v(x_1,x_2,x_3)w(y_1,y_2,y_3).
\]
Hence $\cR(A\otb B\subset M\otb N)\cong \cR\times \cS$, as desired.\hfill$\blacksquare$

\subsection*{Proof of Theorem \ref{F}} 
Let $\mathcal R_1,\mathcal R_2$ be countable ergodic p.m.p. equivalence relations on probability spaces $(X_1,\mu_1), (X_2,\mu_2)$. Let $\mathcal T_1,\mathcal T_2$ be hyperfinite ergodic p.m.p. equivalence relations on probability spaces $(Y_1,\nu_1), (Y_2,\nu_2)$. We assume that $\mathcal R_1$ is strongly ergodic and that $\mathcal R:=\mathcal R_1\times\mathcal T_1\cong\mathcal R_2\times\mathcal T_2$. 
We identify $X_1\times Y_1=X_2\times Y_2$, and denote $M=L(\cR)$, $A=L^\infty(X_1\times Y_1)=L^\infty(X_2\times Y_2)$, $A_1=L^\infty(X_1)$, $A_2=L^\infty(X_2)$, $B_1=L^\infty(Y_1)$, and $B_2=L^\infty(Y_2)$. Lifting the isomorphism between $\cR_1\times \cT_1$ and $\cR_2\times \cT_2$ to the corresponding von Neumann algebras, we get an identification 
\[
[A_1\otb B_1\subset L(\cR_1)\otb L(\cT_1)] = [A_2\otb B_2 \subset L(\cR_2)\otb L(\cT_2)].
\]
Since $\cT_2$ is a hyperfinite ergodic p.m.p. equivalence relation on $(Y_2,\nu_2)$, we can identify $(Y_2,\nu_2)=(Y_0,\nu_0)^{\mathbb N}$, where $Y_0=\{0,1\}$ and $\nu_0=\frac{1}{2}(\delta_0+\delta_1)$, in such a way that $(y_k)\mathcal T_2(z_k)$ if and only if $y_k=z_k$ for all but finitely many $k\in\mathbb N$. We can then identify $B_2=\bar{\otimes}_{k\in\mathbb N}L^{\infty}(Y_0,\nu_0)_k$, and we put $B_{2,n}=\bar{\otimes}_{k\geq n}L^{\infty}(Y_0,\nu_0)_k$, for $n\in\mathbb N$. As in the proof of  Lemma \ref{aux}, we observe that $\prod_\omega B_{2,n}\subset M'\cap A^\omega$. Moreover, by strong ergodicity of $\cR_1$, we also know that $M'\cap A^\omega\subset B_1^\omega$. Putting these inclusions together, we deduce from Lemma~\ref{easy} that $B_{2,n}\subset_{\eps_n} B_1$ where $\eps_n\rightarrow 0$. Since $\eps$-containment does not change when passing to a common superalgebra, and $B_{2,n},B_1\subset A$, choosing $k$ such that $\eps_k<1$ and applying Lemma \ref{epsilon} we get that
\[
B_{2,k} \emb_A B_1.
\]
Hence we can find nonzero projections $p_2\in B_{2,k}$, $p_1\in B_1$, a $^*$-homomorphism $\psi:B_{2,k}p_2\rightarrow B_1p_1$, and a nonzero partial isometry $v\in p_1Ap_2$ such that $\psi(x)v=vx$ for all $x\in B_{2,k}p_2$. Observing that $A$ is abelian and multiplying with $v^*$ on both sides, we find a nonzero projection $p'\in A$ such that
\[
B_{2,k}p' \subset B_1p'.
\]
Now, $B_{2,k}\subset B_2$ is the subalgebra of complex-valued functions on $Y_2 = Y_0^\N$ that do not depend on the first $k-1$ coordinates. Denoting for a fixed $b\in \prod_{i=1}^{k-1} Y_0$ by $e_b:=\mathbbm{1}_{\{b\}\times\prod_{i\geq k} Y_0}$ the indicator function of ${\{b\}\times\prod_{i\geq k} Y_0}$, we thus see that for any such $b$ we have
\[
B_2e_b = B_{2,k}e_b.
\]
As $\sum_{b\in \prod_{i=1}^{k-1} Y_0} e_b = 1_{Y_2}$, we can find $b$ such that $p:=e_b p'\neq 0$. From the above discussion it follows that
\[
B_2p\subset B_1p.
\]
Taking relative commutants, this implies that
\begin{equation}\label{eq:psub}
p(L(\cR_1)\otb B_1)p\subset p(L(\cR_2)\otb B_2)p.
\end{equation}
Since $p\in A$, we can write $p=\mathbbm{1}_Z$ for some measurable subset $Z\subset X_1\times Y_1$ of positive measure. Writing $Z_y=\{x\in X_1\mid (x,y)\in Z\}\subset X_1$ and passing to a subprojection if necessary, we can assume that for every $y\in Y_1$ either $\mu_1(Z_y)=c>0$ or $\mu_1(Z_y)=0$. Choosing a subset $C_1\subset X_1$ with $\mu_1(C_1)=c$, we consider
\[
S:=\{(y,\varphi)\in Y_1\times [\cR_1]\mid \varphi(Z_y)=C_1\}.
\]
Denoting by $\pi_1$ the projection onto the first component, it follows from the ergodicity of $\mathcal R_1$ that $\pi_1(S)=\{y\in Y_1\mid \mu_1(Z_y)\neq 0\}$. Hence by putting $\psi_y=\text{id}_{X_1}$ when $\mu_1(Z_y)=0$, it follows from \cite[Theorem~A.16]{Ta01} that we can find a measurable map $\psi:Y_1\rightarrow [\cR_1]$ such that $\psi_y(Z_y)=C_1$ whenever $\mu_1(Z_y)\neq 0$. Identifying $L(\cR_1)\otb B_1=L^\infty(Y_1,L(\cR_1))$ and putting
\[
u=u_\psi:Y_1\rightarrow L(\cR_1):y\mapsto u_{\psi_y},
\]
we see that $u\in L(\cR_1)\otb B_1$, $uAu^*=A$, and moreover there are projections $q\in A_1$, $\tilde q\in B_1$ such that $upu^* = q\otimes \tilde q$. Conjugating \eqref{eq:psub} by $u$, we thus get
\begin{equation}\label{eq:qsub}
qL(\cR_1)q\,\otb B_1\tilde q \subset (q\otimes\tilde q)u(L(\cR_2)\otb B_2)u^*(q\otimes\tilde q).
\end{equation}
Writing $P=\Ad(u)(L(\cR_2))$, $Q=\Ad(u)(L(\cT_2))$, $A_3=\Ad(u)(A_2)$, and $B_3=\Ad(u)(B_2)$, we get $M=P\otb Q$ and, since $u$ normalizes $A$,
\begin{equation}\label{eq:13}
A_1\otb B_1 = A_2\otb B_2 = A_3\otb B_3.
\end{equation}
Taking relative commutants in \eqref{eq:qsub}, we get $B_3(q\otimes\tilde q)\subset q\otimes B_1\tilde q$. In particular we can  find a von Neumann subalgebra $B_4\subset B_1$ such that $B_3(q\otimes\tilde q)=q\otimes B_4\tilde q$. Taking relative commutants again from both points of view, we get
\[
qL(\cR_1)q\,\otb [(B_4\tilde q)'\cap \tilde qL(\cT_1)\tilde q] = (q\otimes \tilde q)(P\otb B_3)(q\otimes \tilde q).
\]
In particular, since $P$ is a factor, we see that the center of the above algebra equals $B_3(q\otimes \tilde q) = q\otimes B_4\tilde q$. Identifying both with $L^\infty(Y)$ for some probability space $(Y,\nu)$ and disintegrating in the above equality we get
\[
\dint_{Y} qL(\cR_1)q\,\otb N_y\,d\nu(y) = \dint_{Y} \bar q_y P \bar q_y \,d\nu(y),
\]
where we decomposed $N:= (B_4\tilde q)'\cap \tilde qL(\cT_1)\tilde q =\dint_{Y} N_y\,d\nu(y)$, and $q\otimes \tilde q=\dint_{Y} \bar q_y\,d\nu(y)\in P\otb B_3$. It follows from \cite[Theorem~IV.8.23]{Ta01} that the above identification splits, i.e., for almost every $y\in Y$ we necessarily have
\[
qL(\cR_1)q\,\otb N_y = \bar q_y P \bar q_y.
\]
Moreover, we see that $B_4\tilde q\subset B_1\tilde q\subset (B_4\tilde q)'\cap \tilde qL(\cT_1)\tilde q$, so we can decompose $B_1\tilde q = \dint_{Y} B_{1,y}\,d\nu(y)\subset \dint_{Y} N_y\,d\nu(y)=N$, where $B_{1,y}\subset N_y$ is a unital inclusion for all $y$. Now $B_1\subset L(\cT_1)$, and hence also $B_1\tilde q\subset \tilde q L(\cT_1)\tilde q$, is a Cartan subalgebra. Since $B_1\tilde q\subset N\subset \tilde q L(\cT_1)\tilde q$, it follows from \cite{Dy63} that also $B_1\tilde q\subset N$ is a Cartan subalgebra. From \cite[Lemma~2.2]{Sp17} we then deduce that $B_{1,y}\subset N_y$ is a Cartan subalgebra for almost every $y$. Furthermore, it follows from \eqref{eq:13} that
\[
\dint_{Y} A_1q\,\otb B_{1,y} \,d\nu(y) = \dint_{Y} A_3 \bar q_y \,d\nu(y).
\]
This identification again splits by \cite[Theorem~IV.8.23]{Ta01}, i.e. for almost every $y\in Y$ we have $A_1q\,\otb B_{1,y}= A_3 \bar q_y$.

From the above discussion it now follows that we have the following identification of inclusions of Cartan subalgebras, for almost every $y\in Y$:
\[
(A_1q\,\otb B_{1,y} \subset qL(\cR_1)q\,\otb N_y) = (A_3 \bar q_y \subset \bar q_y P \bar q_y).
\]
Writing $\cT_y=\cR(B_{1,y}\subset N_y)$, $r=\uptau(q)$, $s=\uptau(\bar q_y)$, we thus get that for almost every $y$:
\[
\cR_2^{s} \cong \cR(A_2^{s}\subset L(\cR_2)^{s}) \cong \cR(A_3^{s}\subset P^{s})\cong \cR(A_1^r\otb B_{1,y}\subset L(\cR_1)^r\otb N_y) \cong \cR_1^r \times \cT_y,
\]
where the last isomorphism follows from Lemma~\ref{lem:split}. Being a subalgebra of $L(\cT_1)\cong R$, $N$ is amenable, and so $N_y$ is an amenable tracial factor for almost every $y$. Choosing $y$ such that the above hold, \cite{Co75} then implies that $N_y$ is either isomorphic to $\mathbb M_n(\C)$ for some $n\in\N$, or to $R$. In the first case, $\cT_y$ is finite and we can find $t>0$ such that $\cR_2\cong\cR_1^t$. In the second case, $\cT_y$ arises as the equivalence relation of a Cartan inclusion $B\subset R$, and as such is a hyperfinite ergodic p.m.p. equivalence relation by \cite{CFW81}. Hence in this case we have $\cR_2\cong \cR_1\times \cT$ where $\cT$ is a hyperfinite ergodic p.m.p. equivalence relation. This finishes the proof of the theorem.
\hfill$\blacksquare$

\section{New characterisations of property Gamma and strong ergodicity}

\subsection*{Proof of Theorem \ref{G}}
(1) $\Rightarrow$ (2) The proof of this implication builds on an argument due to Popa (see the proof of \cite[Proposition 7]{Oz03}) Assume that $M$ has property Gamma and let $\{x_k\}_{k\geq 1}\subset M$ be a $\|.\|_2$-dense sequence. We construct inductively a sequence $\{B_n\}_{n\geq 1}$ of commuting $*$-subalgebras of $M$ and a projection $p_n\in B_n$ such that $B_n\cong \mathbb M_2(\mathbb C)$, $\uptau(p_n)=\frac{1}{2}$, and $\|x_kp_n-p_nx_k\|_2\leq\frac{1}{2^n}$, for every $n\geq 1$ and $1\leq k\leq n$. Indeed, suppose that $B_1,...,B_n$ and $p_1,...,p_n$ are constructed. Define $C_n=B_1\otimes...\otimes B_n\cong\mathbb M_{2^n}(\mathbb C)$. Since $C_n$ is finite dimensional and $M$ has property Gamma, we can find a projection $p_{n+1}\in C_n'\cap M$ such that $\uptau(p_{n+1})=\frac{1}{2}$, and $\|x_kp_{n+1}-p_{n+1}x_k\|_2\leq\frac{1}{2^n}$, for every $1\leq k\leq n+1.$ Since $C_n'\cap M$ is a II$_1$ factor, we can also find a partial isometry $v_{n+1}\in C_n'\cap M$ such that $v_{n+1}^*v_{n+1}=p_{n+1}$ and $v_{n+1}v_{n+1}^*=1-p_{n+1}$. It is now clear that the algebra $B_{n+1}$ generated by $p_{n+1}$ and $v_{n+1}$ has the desired property.

For $N\geq 1$, let $D_N=\mathbb Cp_N\oplus\mathbb C(1-p_N)$ be the algebra generated by $p_N$. Denote $R=\bar{\otimes}_{i\geq 1}B_i$, $A=\bar{\otimes}_{i\geq 1}D_i$, and  $A_n=\bar{\otimes}_{i\geq n}D_i$, for $n\geq 1$. Then $R$ is a hyperfinite II$_1$ factor and $A$ is a Cartan subalgebra. Fix $n\geq k\geq 1$, and define $A_{n,N}=\bar{\otimes}_{N\geq i\geq n}D_i$, for $N\geq n$. Then for every $N\geq n$ we have  $\|E_{D_N'\cap M}(x_k)-x_k\|_2\leq \frac{1}{2^{N-1}}.$ Since $E_{A_{n,N+1}'\cap M}=E_{A_{n,N}'\cap M}\circ E_{D_{N+1}'\cap M}$, we get that 
$\|E_{A_{n,N+1}'\cap M}(x_k)-E_{A_{n,N}'\cap M}(x_k)\|_2\leq\frac{1}{2^N}$, for every $N\geq n$. Since $\|E_{A_{n,n}'\cap M}(x_k)-x_k\|_2\leq\frac{1}{2^{n-1}}$, by combining these inequalities we get  $\|E_{A_{n,N}'\cap M}(x_k)-x_k\|_2\leq \frac{1}{2^{n-1}}+....+\frac{1}{2^{N-1}}<\frac{1}{2^{n-2}}$, for every $N\geq n$. Since $\lim\limits_{N\rightarrow\infty}E_{A_{n,N}'\cap M}(x_k)=E_{A_n'\cap M}(x_k)$,  we conclude that \begin{equation}\label{contain}\|E_{A_n'\cap M}(x_k)-x_k\|_2\leq\frac{1}{2^{n-2}},\;\;\;\text{for every $n\geq k\geq 1$}.\end{equation}
Next, we claim that if $N$ is a tracial von Neumann algebra and $C\subset N$ is a $*$-subalgebra isomorphic to $\mathbb M_{\ell}(\mathbb C)$, for some $\ell\geq 1$, then $N$ is generated by $C$ and $C'\cap N$. To see this, let $p\in C$ be a projection of trace $\frac{1}{\ell}$ and $w_1,...,w_{\ell}\in C$ be partial isometries such that $w_1=p$, $w_iw_i^*=p$, for all $1\leq i\leq\ell$, and $\sum_{i=1}^{\ell}w_i^*w_i=1$.
If $x\in pNp$, then $\sum_{i=1}^{\ell}w_i^*xw_i\in C'\cap N$, hence $x=p(\sum_{i=1}^{\ell}w_i^*xw_i)p\in C\vee (C'\cap N)$.
If $x\in N$, then $x=\sum_{i,j=1}^{\ell}w_i^*(w_ixw_j^*)w_j$ and since $w_ixw_j^*\in pNp$, for all $1\leq i,j\leq\ell$, we get that $x\in C\vee (C'\cap N)$, which proves our claim.

Finally, since $C_{n-1}=B_1\otimes...\otimes B_{n-1}\cong\mathbb M_{2^{n-1}}(\mathbb C)$ is a $*$-subalgebra of $A_n'\cap M$, the claim implies that $A_n'\cap M= C_{n-1}\vee (C_{n-1}'\cap(A_n'\cap M))$, for every $n\geq 1$.  Since $C_{n-1}\subset R$ and $C_{n-1}'\cap (A_n'\cap M)\subset A'\cap M$, we derive that $A_n'\cap M\subset R\vee (A'\cap M)$, for every $n\geq 1$. In combination with \eqref{contain}, it follows that $x_k\in R\vee (A'\cap M)$, for every $k\geq 1$, and hence that $M=R\vee (A'\cap M)$.

(2) $\Rightarrow$ (3) Assume that $R\subset M$  is a hyperfinite subfactor and $A\subset R$ is a Cartan subalgebra such that $M=R\vee (A'\cap M)$. By Lemma \ref{normalizer} we have that $\mathcal R(A\subset M)=\mathcal R(A\subset R)$ and thus $\mathcal R(A\subset M)$ is hyperfinite and ergodic.

(3) $\Rightarrow$ (1) Assume that $A\subset M$ is a regular von Neumann subalgebra such that $\mathcal R:=\mathcal R(A\subset M)$ is ergodic and hyperfinite. Since $\mathcal R$ is hyperfinite,  we can find a sequence $(a_n)_{n\geq 1}\subset\mathcal U(A)$ such that $\uptau(a_n)=0$, for every $n$, and $\|a_n\circ\alpha-a_n\|_2\rightarrow 0$, for every $\alpha\in[\mathcal R].$ If $u\in\mathcal N_M(A)$, let $\alpha_u\in [\mathcal R]$ such that $uau^*=a\circ\alpha_u$, for every $a\in A$.  Thus,  $\|ua_nu^*-a_n\|_2\rightarrow 0$, for every $u\in\mathcal N_M(A)$. Since $A\subset M$ is regular, we deduce that $(a_n)_n\in\mathcal U(M'\cap M^{\omega})$. Since $\uptau_{\omega}((a_n)_n)=0$, it follows that $M$ has property Gamma.
\hfill$\blacksquare$

\subsection*{Proof of Theorem \ref{H}} (1) $\Rightarrow $ (2) Assume that $\mathcal S$ is strongly ergodic. 
Let $\varphi:X\rightarrow\text{Aut}(\mathcal T)$  be a measurable map satisfying $\varphi(x_1)^{-1}\varphi(x_2)\in[\mathcal T]$, for almost every $(x_1,x_2)\in\mathcal S$. For $y\in Y$, define the map $\varphi_y:X\rightarrow Y$ given by $\varphi_y(x)=\varphi(x)(y)$. Then $(\varphi_y(x_1),\varphi_y(x_2))\in\mathcal T$, for almost every $(x_1,x_2)\in\mathcal S$. Since $\mathcal S$ is strongly ergodic and $\mathcal T$ is hyperfinite, by \cite[Theorem A.2.2]{HK05}, 
we get that $\varphi_y(x)$ is contained in a single $\mathcal T$-class, for almost every $x\in X$.  
In other words, for all $y\in Y$, we have that $(\varphi(x_1)(y),\varphi(x_2)(y))\in\mathcal T$, for almost every $(x_1,x_2)\in X$. By Fubini's theorem, we can therefore find $x_2\in X$ such that $(\varphi(x_1)(y),\varphi(x_2)(y))\in\mathcal T$, for almost every $(x_1,y)\in X\times Y$. Hence, $\varphi(x_1)^{-1}\varphi(x_2)\in [\mathcal T]$, for almost every $x_1\in X$, which proves (2).

(2) $\Rightarrow$ (1)   Assume that $\mathcal S$ satisfies (2), and suppose by contradiction that $\mathcal S$ is not strongly ergodic.   Let $Y_0=\{0,1\}$ together with the probability measure $\nu_0=\frac{1}{2}(\delta_0+\delta_1)$, and define $(Y,\nu)=(Y_0,\nu_0)^{\mathbb N}$. Consider the countable ergodic p.m.p. equivalence relation $\mathcal T$ on $(Y,\nu)$ given by $(y_n)\mathcal T(z_n)$ if and only if $y_n=z_n$, for all but finitely many $n\in\mathbb N$. 
Since $\mathcal S$ is not strongly ergodic, \cite[Theorem 2.1]{JS85} provides a factor map $\pi:(X,\mu)\rightarrow (Y,\nu)$ such that $\pi([x]_{\mathcal S})=[\pi(x)]_{\mathcal T}$, for almost every $x\in X$.

Let $\rho_0=\text{id}_{Y_0}$ and $\rho_1:Y_0\rightarrow Y_0$ be given by $\rho_1(0)=1$ and $\rho_1(1)=0$. 
We define $\rho:Y\rightarrow\text{Aut}(\mathcal T)$ by letting $\rho(y)(t)=(\rho_{y_n}(t_n)),$ for all $y=(y_n)\in Y$ and $t=(t_n)\in Y$.
Then $\rho$ satisfies \begin{equation}\label{rho}\rho(y)^{-1}\rho(z)\in [\mathcal T]\; \Longleftrightarrow \; (y,z)\in\mathcal T,\;\;\text{for all $y,z\in Y$.}\end{equation}
Define $\varphi:X\rightarrow\text{Aut}(\mathcal T)$ by letting $\varphi(x)=\rho(\pi(x))$. 
Then for almost every $(x_1,x_2)\in\mathcal S$, we have that $(\pi(x_1),\pi(x_2))\in\mathcal T$, thus $\varphi(x_1)^{-1}\varphi(x_2)=\rho(\pi(x_1))^{-1}\rho(\pi(x_2))\in [\mathcal T]$. On the other hand, \begin{align*}&(\mu\times\mu)(\{(x_1,x_2)\in X\times X\mid\varphi(x_1)^{-1}\varphi(x_2)\in [\mathcal T]\})\\&=(\mu\times\mu)(\{(x_1,x_2)\in X\times X\mid (\pi(x_1),\pi(x_2))\in\mathcal T\})\\&=(\nu\times\nu)(\mathcal T)=0.\end{align*}
This implies that there does not exist $\alpha\in\text{Aut}(\mathcal T)$ such that $\varphi(x)\alpha\in[\mathcal T]$, for almost every $x\in X$, which contradicts the fact that $\mathcal S$ satisfies condition (2).

(1) $\Rightarrow$ (3) Assume that $\mathcal S$ is strongly ergodic. Then $(L^{\infty}(X)^{\omega})^{[\mathcal S]}=\mathbb C1$, and (3) follows by applying Theorem \ref{rigid} to the trivial factor map. 

(3) $\Rightarrow$ (1)  Assume that $\mathcal S$ satisfies (3), and suppose by contradiction that $\mathcal S$ is not strongly ergodic.  Let $(Y,\nu)=(Y_0,\nu_0)^{\mathbb N}$ and $\mathcal T$ be defined as in the proof of (2) $\Rightarrow$ (1).  
Since $\mathcal S$ is not strongly ergodic, \cite[Theorem 2.1]{JS85} provides a factor map $\pi:(X,\mu)\rightarrow (Y,\nu)$ such that $\pi([x]_{\mathcal S})=[\pi(x)]_{\mathcal T}$, for almost every $x\in X$.

Let $R=\bar{\otimes}_{n\in\mathbb N}\mathbb M_2(\mathbb C)$ be the hyperfinite II$_1$ factor. Let $u\in\mathbb M_2(\mathbb C)$ be a unitary with $\uptau(u)=0$.
Let $\rho_0=\text{id}_{\mathbb M_2(\mathbb C)}$ and $\rho_1=\text{Ad}(u)\in\text{Aut}(\mathbb M_2(\mathbb C))$. We define $\rho:Y\rightarrow\text{Aut}(R)$ by letting $$\rho(y)=\otimes_{n\in\mathbb N}\;\rho_{y_n}\;\;\;\;\text{for all $y=(y_n)\in Y$.}$$

If $y,z\in Y$ and $(y,z)\in\mathcal T$, then we clearly have that $\rho(y)^{-1}\rho(z)\in\text{Inn}(R)$. Conversely, assume that $\rho(y)^{-1}\rho(z)\in\text{Inn}(R)$, for some $y,z\in Y$. 
If $(y,z)\notin\mathcal T$, then we can find a subsequence $(k_n)$ on $\mathbb N$ such that $y_{n_k}\not= z_{n_k}$, for all $k\geq 1$. 
Since $u\notin\mathbb C1$ we can find  $v\in\mathbb M_2(\mathbb C)$ such that $\rho_1(v)=uvu^*\not=v$. For $k\geq 1$, let $v_k=\otimes_{n\in\mathbb N}v_{n,k}\in R$, where $v_{n,k}=v$, if $n=n_k$, and $v_{n,k}=1$, if $n\not=n_k$. Then $$\|\rho(y)^{-1}\rho(z)(v_k)-v_k\|_2=\|\rho_{y_{n_k}}\rho_{z_{n_k}}^{-1}(v)-v\|_2=\|uvu^*-v\|_2>0,\;\;\;\text{for all $k\geq 1$.} $$
Since $(v_k)_k\in R'\cap R^{\omega}$ this contradicts the fact that $\rho(y)^{-1}\rho(z)\in\text{Inn}(R)$. Altogether, we have shown that $\rho$ satisfies $\rho(y)^{-1}\rho(z)\in\text{Inn}(R)\; \Longleftrightarrow \; (y,z)\in\mathcal T$, for all $y,z\in Y$.

Finally, define $\theta:X\rightarrow\text{Aut}(R)$ by letting $\theta(x)=\rho(\pi(x))$. By repeating the argument from the end of the proof of (2) $\Rightarrow$ (1) one derives a contradiction with the fact that $\mathcal S$ satisfies condition (3). \hfill$\blacksquare$


\begin{thebibliography}{ABC99}
\bibitem[CD18]{CD18} I. Chifan and S. Das: {\it A remark on the ultrapower algebra of the hyperfinite factor}, preprint arXiv:1802.06628, to appear in Proc. Amer. Math. Soc.
\bibitem[Ch82]{Ch82} M. Choda: {\it Inner amenability and fullness}, Proc. Amer. Math. Soc. {\bf 86} (1982), 663-666.
\bibitem[Co75]{Co75} A. Connes: {\it Classification of injective factors}, Ann. of Math. {\bf 74} (1976),
73-115.
\bibitem[CFW81]{CFW81} A. Connes, J. Feldman and B. Weiss: {\it An amenable equivalence relation is generated by a single transformation}, Ergodic Theory and Dynam. Systems, {\bf 1} (1981), 431-450.
\bibitem[CJ82]{CJ82} A. Connes and V.F.R. Jones: {\it A II$_1$ factor with two nonconjugate Cartan subalgebras}, Bull. Amer. Math. Soc. (N.S.) {\bf 6} (1982), no. 2, 211-212. 
\bibitem[DHI16]{DHI16} D. Drimbe, D. Hoff, and A. Ioana: {\it Prime II$_1$ factors arising from irreducible lattices in products of rank
one simple Lie groups}, preprint arXiv:1611.02209, to appear in J. Reine. Angew. Math.
\bibitem[DV16]{DV16} T.  Deprez and S. Vaes: {\it Inner amenability, property Gamma, McDuff II$_1$ factors and stable equivalence relations}, preprint arXiv:1604.02911, to appear in Ergodic Theory and Dynam. Systems.
\bibitem[Dy63]{Dy63} H.A. Dye: {\it On groups of measure preserving transformations, II}, Amer. J. Math. {\bf 85} (4) (1963) 551-576.
\bibitem[EK97]{EK97} D. Evans and Y. Kawahigashi: {\it Quantum symmetries on operator algebras}, Oxford Mathematical Monographs. Oxford Science Publications. The Clarendon Press, Oxford University Press, New York, 1998. xvi+829 pp.
\bibitem[FM77]{FM77} J. Feldman and C.C. Moore: {\it Ergodic Equivalence Relations, Cohomology, and Von Neumann Algebras. I, II.}
Trans. Amer. Math. Soc. {\bf 234} (1977), 289-324, 325-359.
\bibitem[Ga99]{Ga99} D. Gaboriau: {\it Co\^{u}t des relations d'\'{e}quivalence et des groupes}, Invent. Math., {\bf 139}
(2000), 41-98.
\bibitem[Gl03]{Gl03} E. Glasner: {\it Ergodic theory via joinings}, Amer. Math. Society, 2003.
\bibitem[Ha79]{Ha79} U. Haagerup: {\it An example of a non-nuclear C$^*$-algebra, which has the metric approximation property}, Invent.
Math. {\bf 50} (1978/79), 279-293.
\bibitem[HK05]{HK05} G. Hjorth and A. Kechris: {\it Rigidity theorems for actions of product groups and countable Borel equivalence
relations}, Mem. Amer. Math. Soc. 177 (2005), no. 833, viii+109 pp.
\bibitem[HMV16]{HMV16} C. Houdayer, A. Marrakchi, and P. Verraedt: {\it Fullness and Connes' $\tau$ invariant of type III tensor product factors}, preprint arXiv:1611.07914, to appear in J. Math. Pures Appl.
\bibitem[Ho15]{Ho15} D. Hoff: {\it Von Neumann Algebras of Equivalence Relations with Nontrivial One-Cohomology}, J. Funct. Anal. {\bf 270} (2016),  1501-1536.
\bibitem[JS85]{JS85} V.F.R. Jones and K. Schmidt: {\it Asymptotically invariant sequences and approximate finiteness}, 
Amer. J. Math. {\bf 109} (1987), no. 1, 91-114.
\bibitem[Ki12]{Ki12} Y. Kida: {\it Inner amenable groups having no stable action}, Geom. Dedicata {\bf 173} (2014), 185-192.
\bibitem[Ki14]{Ki14} Y. Kida: {\it Splitting in orbit equivalence, treeable groups, and the Haagerup property},
 Hyperbolic geometry and geometric group theory, 167-214, Adv. Stud. Pure Math., {\bf 73}, Math. Soc. Japan, Tokyo, 2017.
 \bibitem[Ki16]{Ki16} Y. Kida: {\it Stable actions and central extensions}, Math. Ann. {\bf 369} (2017), no. 1-2, 705-722.
 \bibitem[Ma17]{Ma17} A. Marrakchi: {\it Stability of products of equivalence relations}, preprint arXiv:1709.00357, to appear in Compos. Math.
 \bibitem[Mc69]{Mc69} D. McDuff: {\it Central sequences and the hyperfinite factor}, Proc. London Math. Soc. {\bf 21} (1970), 443-461.
 \bibitem[MvN36]{MvN36} F. Murray and J. von Neumann: {\it On rings of operators}, Ann. of Math., {\bf 37} (1936),
116-229.
 \bibitem[MvN43]{MvN43} F. Murray and J. von Neumann: {\it Rings of operators}, IV. Ann. of Math. {\bf 44} (1943), 716-808.
\bibitem[OW80]{OW80} D. Ornstein and B. Weiss: \textit{Ergodic theory of amenable group actions. I. The Rohlin lemma.} Bull. Amer. Math. Soc. (N.S.), 2(1):161-164, 1980.
\bibitem[Oz03]{Oz03} N. Ozawa: {\it Solid von Neumann algebras}, Acta Math. {\bf 192} (2004), no. 1, 111-117.
\bibitem[Po01]{Po01} S. Popa: {\it On a class of type II$_1$ factors with Betti numbers invariants}, Ann. of Math. {\bf 163} (2006), 809-899.
 \bibitem[Po03]{Po03} S. Popa: {\it Strong rigidity of II$_1$ factors arising from malleable actions of w-rigid groups I}, Invent. Math. {\bf 165} (2006), 369-408.
 \bibitem[Po04]{Po04} S. Popa: {\it Deformation and rigidity in the study of II$_1$ factors}, Mini-Course at College
de France, November 2004.
  \bibitem[Po06a]{Po06a} S. Popa: {\it On the superrigidity of malleable actions with spectral gap}, J. Amer. Math. Soc. {\bf 21} (2008), 981-1000.
 \bibitem[Po06b]{Po06b} S. Popa: {\it On Ozawa's property for free group factors}, Int. Math. Res. Not. IMRN 2007, no. 11, Art.
ID rnm036, 10 pp.
\bibitem[Po09]{Po09} S. Popa: {\it On spectral gap rigidity and Connes invariant $\chi(M)$}, Proc. Amer. Math. Soc. {\bf 138} (2010), no. 10, 3531-3539.
 \bibitem[Po13]{Po13} S. Popa: {\it Independence properties in subalgebras of ultraproduct II$_1$ factors}, J. Funct. Anal. {\bf 266} (2014), 5818-5846.
 \bibitem[PS03]{PS03} S. Popa and D. Shlyakhtenko: {\it Cartan subalgebras and bimodule decompositions of II$_1$ factors}, Math. Scand. {\bf 92} (2003), 93-102.
\bibitem[PV11]{PV11} S. Popa and S. Vaes: {\it Unique Cartan decomposition for II$_1$ factors arising from arbitrary actions of free groups}, Acta Math. {\bf 212} (2014), 141-198.
\bibitem[Sc85]{Sc85} K. Schmidt: {\it Cohomology and the absence of strong ergodicity for ergodic group actions}. Operator algebras and their connections with topology and ergodic theory (Bu\c{s}teni, 1983), 486-496, Lecture Notes in Math., 1132, Springer, Berlin, 1985.
\bibitem[Sp17]{Sp17} P. Spaas: {\it Non-classification of Cartan subalgebras for a class of von Neumann algebras}, Adv. in Math. {\bf 332} (2018), 510-552.
 \bibitem[Ta01]{Ta01} M. Takesaki: {\it Theory of Operator Algebras I}, ser. Encyclopaedia of Mathematical Sciences. Springer, {\bf 124}, 2001. xix+415 pp.
 \bibitem[TD14]{TD14} R. Tucker-Drob: {\it Invariant means and the structure of inner amenable groups}, preprint arXiv:1407.7474.
 \bibitem[Va08]{Va08} S. Vaes: {\it Explicit computations of all finite index bimodules for a family of II$_1$ factors}, Ann. Sci. \'{E}c. Norm.
Sup. {\bf 41} (2008), 743-788.
\end{thebibliography}
\end{document}